%% file: main.tex
\documentclass{cocv}

\input{style.sty}

\usepackage{xstring}
\usepackage[normalem]{ulem}
\definecolor{darkgreen}{rgb}{0.0, 0.5, 0.0}
\newcommand{\edit}[3][1]{%
	\IfEq{#1}{2}{\def\mysecondvar{}}{\def\mysecondvar{#2}}
	\ifx\mysecondvar\empty{}\else{%
		\IfEq{#1}{1}{
			\ifmmode%
			\text{\color{red}\sout{\ensuremath{#2}}}%
			\else%
			{\color{red} \sout{#2}}
			\fi%
		}{\color{red}#2}\ %
	}\fi{\color{darkgreen}#3}}
\renewcommand{\edit}[3][1]{#3}

\begin{document}
	\title{Dynamical Programming for off-the-grid dynamic Inverse Problems}
	\thanks{This work was funded by the ANR CIPRESSI project, grant ANR-19-CE48-0017-01 of the French Agence Nationale de la Recherche.}
	
	\author{Vincent Duval}
	\address{INRIA-Paris, MOKAPLAN, 75012 Paris, France (\email{vincent.duval@inria.fr\ \&\ robert.tovey@inria.fr})
		\\
		CEREMADE, CNRS, UMR 7534, Universit\'e Paris-Dauphine, PSL University 75016 Paris, France}
	\author{Robert Tovey}\sameaddress{1}
	
	\begin{abstract}
		In this work we consider off-the-grid algorithms for the reconstruction of sparse measures from time-varying data. In particular, the reconstruction is a finite collection of Dirac measures whose locations and masses vary continuously in time. Recent work showed that this decomposition was possible by minimising a convex variational model which combined a quadratic data fidelity with dynamical Optimal Transport. We generalise this framework and propose new numerical methods which leverage efficient classical algorithms for computing shortest paths on directed acyclic graphs. Our theoretical analysis confirms that these methods converge to globally optimal reconstructions. Numerically, we show new examples for unbalanced Optimal Transport penalties, and for balanced examples we are 100 times faster in comparison to the previously known method.
	\end{abstract}
	
	\begin{resume}
		Dans cet article, nous considérons des algorithmes de reconstruction sans grille pour des mesures parcimonieuses à partir de données dynamiques. En particulier, le signal reconstruit est une somme finie de mesures de Dirac dont les positions et les amplitudes varient continûment.
		Des travaux récents ont montré qu'une telle décomposition était possible dans un cadre variationnel convexe qui combine une attache aux données quadratique avec la formulation dynamique du transport optimal. Nous généralisons ce cadre et proposons une nouvelle méthode numérique qui exploite des algorithmes efficaces classiques de plus court chemin sur des graphes acycliques orientés. Notre analyse théorique confirme que ces méthodes convergent vers des minimiseurs globaux. Nous illustrons ces méthodes sur des exemples nouveaux impliquant le transport optimal non-équilibré, et, pour le transport classique, notre méthode est près de 100 fois plus rapide que l'état de l'art.
	\end{resume}

	
	\maketitle
	\layout{\newpage}
		
	\IfFileExists{sections/intro.tex}{
		\input{sections/intro.tex}
		\input{sections/preliminary.tex}
		\input{sections/generalisation.tex}	
		\input{sections/fwcvter.tex}
		\input{sections/generaldp.tex}
		\input{sections/WFRexample.tex}	
		\input{sections/results.tex}
		\input{sections/conclusion.tex}
	}{

\input{intro.tex}\input{preliminary.tex}\input{generalisation.tex}\input{fwcvter.tex}\input{generaldp.tex}\input{WFRexample.tex}\input{results.tex}\input{conclusion.tex}
	}
	
	\bibliographystyle{abbrv}
	\bibliography{bib}
	
	\appendix
	\IfFileExists{sections/intro.tex}{
		\input{appendices/flat_metric.tex}
		\input{appendices/Theta_thm.tex}
		\input{appendices/lemmas.tex}
		\input{appendices/energy_results.tex}
	}{

\input{flat_metric.tex}\input{Theta_thm.tex}\input{lemmas.tex}\input{energy_results.tex}
	}
	
\end{document}

%% file: intro.tex

\section{Introduction}

The signal-processing task of dynamical super-resolution involves retrieving fine-scale features, in space and time, from a signal which evolves over time. A convex variational model was recently proposed for such tasks using Optimal Transport (OT) to regularise the associated inverse problem \citep{BrediesBalanced, BrediesOptimal}. This new approach allows the decomposition of a signal into a finite sum of smooth curves, for example to track the centers of multiple particles in time with smooth trajectories. Similar ideas were explored for a specific example in \cite{Alberti2019} where the shape of curves is built into the model, and without appealing to OT.

In this work we focus on dynamical super-resolution problems regularised by OT. We can consider potential models to be partitioned into two classes, depending on whether the particles have constant mass/brightness in time, or if mass is allowed to vary. These classes are referred to as balanced or unbalanced problems respectively, mathematically encoded in the choice of OT cost. Current literature provides analysis for the balanced Benamou-Brenier (BB) \citep{BrediesBalanced} and the unbalanced Wasserstein-Fisher-Rao (WFR) \citep{BrediesUnbalanced} energies, both are shown to reconstruct data into a finite number of smooth curves with constant or smoothly varying mass. Initial numerical experiments for the Benamou-Brenier model have also been carried out showing great promise \citep{BrediesNumerics}, however current methods are too slow for large-scale applications.

Let $\Omega$ be an open bounded spatial domain. At the heart of the analysis of Bredies \etal\ is the interplay between measures $\dm(t,\x)$ defined on the time-space cylinder $\ci01\times\bar\Omega$, and measures $\cm(\hg)$ defined on the space of curves $\hg = (\h,\g)$ with mass $\h\in\Co(\ci01)$ and trajectory $\g\in\Co(\ci01;\bar\Omega)$. We can think of $\dm$ as representing the evolving physical volume that can be observed with a microscope, whereas $\cm$ more efficiently represents a collection of (trajectories of) particles, that we wish to reconstruct. 
The structure of the problem as seen from this second viewpoint closely resembles the Beurling-LASSO which is now well-understood \citep{brediesInverseProblemsSpaces2013,azais_spike_2015,duval-focm-2015,poon_geometry_2020}. Its main advantage is that it paves the way for ``off-the-grid'' numerical methods when solving such dynamical inverse problems. Recent works in the field of sparse spike recovery \citep{brediesInverseProblemsSpaces2013,boydAlternatingDescentConditional2017,denoyelleSlidingFrankWolfeAlgorithm2019} have demonstrated that it is possible to design efficient numerical solvers without reconstructing the unknown on a grid, by exploiting a conditional gradient descent / Frank-Wolfe approach together with a good knowledge of the regularising term. Indeed, the Frank-Wolfe minimisation algorithm and its variants (see the review \cite{Jaggi2013}) build iterates that are convex combinations of the extreme points of level sets of the regulariser; being able to easily encode and handle such extreme points makes it possible to solve variational problems in a continuous (or up to floating point) setting.
Moreover, having iterates that are convex combinations of a few extreme points of the level sets of the regulariser is particularly relevant, as it is known in inverse problems that some solutions have precisely that structure when the observation consists of a finite number of linear measurements \citep{unser_splines_2017,bredies_sparsity_2019,Boyer2019}. 

\edit{}{
\subsection{Motivating example}\label{sec:BB example intro}
To make these observations and the contribution of this work more concrete we will make reference to the motivating example described in \cite{BrediesBalanced} with numerical examples in \cite{BrediesNumerics}. We will describe this problem briefly here and postpone precise assumptions and details of function spaces to Section~\Rref{sec:WFR example} where it is the $\delta=+\infty$ case. 
For $\alpha,\beta>0$ the Benamou-Brenier penalty $\creg$ is a map of non-negative space-time measures $\dm$ defined by
\begin{equation}\label{eq: BB energy}
	\creg(\dm) \eqdef \inf_{v\in \xLtwo_\dm(\ci01\times\bar\Omega;\xR^d)}\Bigg\{\int_\X \left[\alpha + \frac\beta2|v|^2\right]\xdif\dm \ \st\ \partial_t\dm + \op{div}(v\dm)=0\Bigg\},
\end{equation}
where the continuity equation is satisfied in the weak sense which will be clarified in \eqref{eq: continuity equation}. The function with $\alpha=0$, $\beta=1$ is referred to as the Benamou-Brenier energy whose main properties can be found in~\cite[Sec. 5.3.1]{santambrogio_optimal_2015}. It was shown (\cf \cite{santambrogio_optimal_2015,BrediesBalanced}) that whenever $\creg(\dm)<+\infty$, there exists a ``disintegration'' into spatial measures $\{\dm_t\}_{t\in\ci01}$ such that the curve $t\mapsto \dm_t$ is continuous in the narrow topology and such that
\begin{equation}
	\forall \psi \in \xLone(\ci{0}{1}\times\bar{\Omega};\dm),\quad \int_\X\psi(t,x) \xdif \dm(t,x) = \int_{0}^1 \left(\int_{\bar{\Omega}}\psi(t,x)\xdif \dm_t(x) \right)\xdif t.
\end{equation}
Thanks to this property, given times $0=t_0<t_1<\ldots<t_T=1$, it is justified to consider ``slices'' $\dm_{t_j}$ and to assume that we are given data $\data_j\in\xR^\m$ (for simplicity assume $\m$ does not change with $j$) at time $t_j$, corresponding to linear observations of $\dm_{t_j}$, given by narrowly continuous linear operators $A_j$. The more involved case of continuous-time observations is handled in \cite{BrediesOptimal} although we do not discuss it further in this work. To solve the corresponding inverse problem, \cite[(43)]{BrediesBalanced} proposes the minimisation of
\begin{equation}\label{bbexample-eq-energy}
	\C{E}(\dm) = \frac12\sum_{j=0}^T\norm{A_j\dm_{t_j}-\data_j}_2^2 + \creg(\dm).
\end{equation}
It is then shown in \cite[Theorem 10]{BrediesBalanced} that there is a minimiser $\dm^*$ which is a finite sum of extreme points of the Benamou-Brenier unit ball, \ie
\begin{equation}\label{bbexample-eq-structuremin-rho}
	\text{for some }a_i\geq 0,\ \g^i\in \ACtwo(\ci01;\bar\Omega),\qquad \forall t \in \ci{0}{1},\qquad \dm^*_t = \sum_{i=1}^{\m(T+1)} a_i \delta_{\g^i(t)},
\end{equation}
where $\ACtwo$ is the set of absolutely continuous function such that their (a.e. defined) pointwise derivative is square-integrable (see also \cite[Sec. 1.1]{Ambrosio2008}).
The sparse structure of $\dm^*$ motivates us to look at problems like \eqref{bbexample-eq-energy} from the viewpoint of measures on paths, \ie\ to have a minimiser $\cm^* = \sum_{i=1}^{\m(T+1)} a_i\delta_{\g^i}$. This will enable us to greatly generalise the types of regulariser $\creg$ for which minimisers are sparsely supported on curves, and develop a new specialised tool for ``off-the-grid'' numerical methods in this setting.
}

\subsection{Outline of the paper}
In Section~\Rref{sec: preliminary} we gather and present several results on the measure representation of the solutions of the continuity equation that are useful to understand the perspective of measures on paths.
Our first main theoretical contribution is found in Section~\Rref{sec: generalisation} where we introduce our main variational inverse problem and provide several theoretical results on the existence of sparse minimisers and their ``geodesic'' structure. Algorithmic details are discussed in Section~\Rref{sec: FW convergence} where we analyse a new stochastic Frank-Wolfe variant which is practically computable and proven to converge almost surely. This analysis is our other main theoretical contribution. 

Numerical work on the Benamou-Brenier example was recently presented in \cite{BrediesNumerics}. The most time consuming step is to compute new extreme points, \ie\ curves $\hg$, to add to the reconstruction. Even though only a small subset of $\xLinfty(\ci01;\xR\times\bar\Omega)$ is involved, optimising over such a set is a challenging task, the computation times are hardly compatible with practical applications. 
Our main numerical contribution is to reformulate this step into a shortest-path problem on an ordered, weighted, directed, acyclic graph in Section~\Rref{sec: dynamical programming}, which can be solved very efficiently using dynamical programming. Numerical experiments demonstrate the efficiency of this new technique in Section~\Rref{sec: results}. In particular, we compare with the algorithm of \cite{BrediesNumerics} for the Benamou-Brenier penalty outlined in Section~\Rref{sec:BB example intro}. Finally we show new results for the unbalanced Wasserstein-Fisher-Rao penalty which is described in detail in Section~\Rref{sec:WFR example}.

\subsection{Notation}
\Paragraph{Convex sets and extreme points.}
Let $V$ be a linear space. For all $\cm_0, \cm_1 \in V$, we define the closed line segment between $\cm_0$ and $\cm_1$ as $\ci{\cm_0}{\cm_1}\eqdef \enscond{\lambda \cm_0+(1-\lambda)\cm_1}{0\leq \lambda\leq 1}$. Similarly, we define the open line segment $\oi{\cm_0}{\cm_1}\eqdef \ci{\cm_0}{\cm_1}\setminus \{\cm_0,\cm_1\}$.
A set $D \subseteq V$ is called \emph{convex} if $\oi{\cm_0}{\cm_1}\subset D$ for all $\cm_0,\cm_1\in D$. We say that $\cm \in  D$ is an \emph{extreme point} (or \emph{atom}) of $D$, and write $\cm\in\Ext D$,  if there are no points $\cm_0, \cm_1 \in D$ such that $\cm \in \oi{\cm_0}{\cm_1}$.
	In other words,
	\begin{equation}\label{eq: extreme point equation}
		\forall \lambda \in \oi01, \forall \cm_0,\cm_1\in D,\quad 		\left(\cm = \lambda \cm_0 + (1-\lambda)\cm_1\right) \quad\implies\quad \left(\cm_0=\cm_1=\cm\right).
	\end{equation}

	Furthermore, it is possible to define the notion of face (and elementary face) of a convex set, which extends the notion of extreme point to higher-dimensional sets. We refer to~\cite{dubins_extreme_1962} for more detail.
	

\Paragraph{Measure spaces.}
For a separable metric space $\HG$ and Banach space $X$, we define $\Cb(\HG;X)$ to be the set of continuous bounded functions from $\HG$ to $X$. When $X=\xR$, we simply write $\Cb(\HG)$. Recall that for any Borel measure $\cm$ on $\HG$, we can define the  non-negative Borel measure $\abs\cm\in\MC^+(\HG)$ by
\begin{equation}
	\abs\cm(A) \eqdef \sup \enscond{\sum_{i=1}^{n} \abs{\cm(A_i)}}{n \in \xN, \{A_1, \ldots, A_n\}\ \mbox{Borel partition of $A$}}
\end{equation}
for all Borel measurable sets $A\subset\HG$. We denote by $\MC(\HG)$ the space of signed Borel measures $\cm$ with finite total variation, \ie\ $\norm\cm\eqdef \int_\HG\xdif\abs\cm <+\infty$.
The total variation $\norm\cdot$ defines a norm on $\MC(\HG)$, but it is sometimes more convenient to use the \emph{narrow topology}, \ie\ the weakest topology on $\MC(\HG)$ which makes the integration against continuous bounded functions a continuous linear form. The narrow topology is equivalent to the weak-* topology on $(\Cb(\HG))'$, in particular, a sequence $\{\cm^n\}_{n\in\xN}\subset\MC(\HG)$ converges to $\cm^*\in\MC(\HG)$ in the narrow topology (denoted $\cm^n\stackrel{*}{\rightharpoonup}\cm$) if
\begin{equation}
	\forall \phi\in \Cb(\HG),\quad
\lim_{n \to +\infty} \int_{\HG} \phi\xdif\cm^n = \int_{\HG} \phi\xdif\cm^*.
\end{equation}
The support of $\cm\in\MC^+(\HG)$ is defined as
\begin{equation}
	\op{supp}(\cm) \eqdef \left(\bigcup\enscond{U}{\cm(U)=0,\ U\text{ is open}}\right)^c.
\end{equation}
This is a closed set satisfying $\cm(\op{supp}(\cm)) = \cm(\HG)$.

\Paragraph{Function domains.}
In this work we consider two measure domains, the time-space cylinder
\begin{equation}
	\X \qtext{for an open, bounded, convex domain }\Omega \subseteq \xR^\d,\ d\geq1
\end{equation}
and a closed set $\HG$ of continuous (weighted) curves which are viewed as pairs $\hg=(\h,\g)$ where $\h(t)$ is the mass at time $t$ and $\g(t)$ is the location. A formal definition will be given in Lemma~\Rref{lmma: metric properties}. In the time-space cylinder we will denote measures $\dm\in\MC(\X)$ with test functions $\psi\in\Co(\X)$, and similarly on the space of curves $\cm\in\MC(\HG)$ and $\phi\in\Cb(\HG)$. 

\Paragraph{Narrowly continuous measures.}
An important subspace of $\MC(\X)$ is the space of narrowly continuous measures. With a slight abuse of the standard notation, we will say $\dm\in\CwX\subset \MC(\X)$ if there exists a map $t\mapsto \dm_t\in \MC(\bar\Omega)$ (informally, $\dm_t = \dm(t,\cdot)$ is a ``time slice'' of $\dm$) such that
\begin{equation}
	\forall \psi\in\Co(\bar\Omega),\qquad \left[t\mapsto \int_{\bar\Omega}\psi(\x)\xdif\dm_t(\x)\right]\in \Co(\ci01)
\end{equation}
and
\begin{equation}\label{eq: disintegration equation}
	\forall \psi\in \xLone_\dm(\X),\qquad \int_\X\psi(t,\x)\xdif\dm(t,\x) = \int_0^1\left(\int_{\bar\Omega}\psi(t,\x)\xdif\dm_t(\x)\right)\xdif t.
\end{equation}
Given a measure on paths $\cm\in\MC(\HG)$ such that $\int_{\HG}\norm{\h}_\infty\xdif \sigma(\h,\g)<+\infty$, one may define the family of measures $\cmt{t}\in\MC(\bar\Omega)$, for $t \in \ci01$, by 
	\begin{equation}\label{eq:pushfwdeval}
	\forall \psi\in\Co(\bar\Omega),\qquad \int_{\bar\Omega}\psi(\x)\xdif [\cmt{t}](\x) \eqdef \int_\HG \h(t)\psi(\g(t))\xdif\cm(\h,\g).
\end{equation}
Formally, $\cmt{t}$ is the image measure of $\cm$ by the evaluation at time $t$. That family is narrowly continuous, and as we explain in Theorem~\Rref{thm: Theta properties} below, it is the evaluation (disintegration) of some measure $\Theta(\sigma) \in \CwX$.

\Paragraph{The continuity equation.} In the rest of the paper, we use the following distributional definition of the continuity equation, formally
\begin{equation}
	\partial_t\dm +\op{div}(\dm v) = g\dm,
\end{equation}
which expresses mass variation (or mass conservation if $g= 0$).
\begin{dfntn}
	Let $\dm\in\MC(\X)$ be a measure. We say that $\dm$ satisfies the \emph{continuity equation} if there exists $v\in \xLone_{\abs{\dm}}(\X;\xR^d)$, $g\in \xLone_{\abs{\dm}}(\X)$ such that
	\begin{equation}\label{eq: continuity equation}
		\forall\psi\in\Co_c^1(\oi01\times\bar\Omega),\qquad \int [\partial_t\psi+\nabla\psi\cdot v + \psi g]\xdif\dm = 0.
	\end{equation}
\end{dfntn}


%% file: preliminary.tex

\section{Preliminary results}\label{sec: preliminary}
As previously mentioned, the main function space of this work is the space of measures on paths. In particular, $\MC(\HG)$ where $\HG\subset \HGF$ is the set of continuous weighted paths in $\bar\Omega$, modelling particles with (varying) mass $\h(t)$ at location $\g(t)$ defined by
\begin{equation}
	\HGF=\enscond{\hg=(\h,\g)}{\h\in\Co(\ci01),\ \g \colon\ci01\rightarrow\bar\Omega,\ \g|_{\{\h\neq0\}} \text{ is continuous}}.
\end{equation}
For technical reasons we permit curves $\g$ which may not be continuous at points $t$ where $\h(t)=0$. Intuitively, the location $\g(t)$ is not necessarily meaningful if the particle has no mass and cannot be observed. 

In this section we review the necessary assumptions for $\MC(\HG)$ to be a sufficiently well-behaved space. Firstly we require $\HG$ to be a complete separable metric space. \edit{}{We follow the suggestion of \cite[Prop. 3.6]{BrediesUnbalanced} who used the flat metric on the space of measures $\enscond{\h\delta_\g\in\MC^+(\X)}{(\h,\g)\in\HG,\ \h\geq 0}$ for a particular $\HG\subset\HGF$. We use the isometric space $(\HGF,\metric)$, the properties of which are given by the following lemma.}
\begin{lmm}\label{lmma: metric properties}
	Define $\metric\colon\HGF\times\HGF\to\hc{0}{+\infty}$ by
	\begin{align}
		\metric((\h_1,\g_1),(\h_2,\g_2)) &\eqdef \sup_{t\in\ci01} \op{d}_F((\h_1(t),\g_1(t)),(\h_2(t),\g_2(t))) \qquad\text{where}\\
		\op{d}_F((r_1,\x_1),(r_2,\x_2)) &\eqdef \splitln{|r_1|+|r_2|}{r_1r_2\leq 0\text{ or } |\x_1-\x_2|\geq 2}{|r_1-r_2| + \min(|r_1|,|r_2|)|\x_1-\x_2|}{\qquad\text{else,}}
	\end{align}
	then $(\sfrac\HGF\sim,\metric)$ is a complete separable metric space where 
	\begin{equation}
		(\h_1,\g_1)\sim (\h_2,\g_2) \iff \h_1=\h_2\qtext{and}\forall t\in\{\h_1\neq0\},\quad \g_1(t)=\g_2(t).
	\end{equation}
	Convergence of a sequence $\hg_n=(\h_n,\g_n)\in\HGF$ in the metric $\metric$ can equivalently be stated as:
	\begin{equation}
		\left[\hg_n\stackrel{\metric}{\to}(\h,\g)\right] \iff \Big[\h_n\to\h \text{ in } \Co(\ci01)\text{ and for all $\epsilon>0$, } \g_n\to\g \text{ in } \Co(\{|\h|\geq\epsilon\})\Big]
	\end{equation}
	Furthermore, for any $\psi\in\Co(\X)$, we have $\Psi\in\Co(\ci01\times\HGF)$ where
	\begin{equation}
		\forall t\in\ci01,\ (\h,\g)\in \HGF, \qquad \Psi(t,\h,\g)\eqdef \h(t)\psi(t,\g(t)).
	\end{equation}
\end{lmm}
\edit{}{The proof is elementary but given in Appendix~\Rref{app: prelims} for completeness as no specific reference could be found.}

A key analytical tool in related prior works (\cf \cite{BrediesBalanced,BrediesOptimal,BrediesUnbalanced,BrediesNumerics}) is a mapping between measures $\dm\in\MC(\X)$ which satisfy the continuity equation, and measures on (weighted) paths $\cm\in\MC(\HGF)$, making some structures become more apparent. We will now recap and expand on those previous results. \edit{}{The first theorem collects results from \cite{Ambrosio2008,BrediesUnbalanced} and provides minor extensions for the scope of the current work. In particular, we remove the necessity for $\h\geq0$ or elements $(\h,\g)$ to satisfy further smoothness conditions.}

\begin{thrm}\label{thm: Theta properties}
	Let $\cm\in\MC(\HGF)$. If $\int_\HGF\norm\h_1\xdif|\cm|(\h,\g)<+\infty$, then there is a unique finite Borel measure $\Theta(\cm)\in\MC(\X)$ such that
	\begin{equation}\label{eq: Theta def.}
		\forall\psi\in\Co(\X),\qquad \int_\X\psi(t,\x)\xdif\Theta(\cm)(t,\x) = \int_\HGF \left(\int_0^1\h(t)\psi(t,\g(t))\xdif t\right)\xdif\cm(\h,\g).
	\end{equation}
	Moreover, 
	\begin{enumerate}
		\item The mapping $\Theta\colon \enscond{\cm\in\MC(\HGF)}{\int_\HGF\norm\h_1\xdif|\cm|<+\infty} \to\MC(\X)$ is linear.
		\item Equality~\eqref{eq: Theta def.} holds for all $\psi\in \xLone_{|\Theta(\cm)|}(\X)$.
		\item If $\int_\HGF\norm\h_\infty\xdif\abs{\cm}<+\infty$, then $\Theta(\cm)\in\CwX$.
		\item Suppose $\h,\g\in \ACtwo(\ci01)$ for $\cm$-a.e. $(\h,\g)\in\HGF$. If there exist Borel measurable functions $v\colon \X\rightarrow \xR^d$ and $g:\X\rightarrow \xR$ such that 
		\begin{align}
			&\h'(t) = g(t,\g(t))\h(t) \text{ for $\cm$-a.e. $(\h,\g)$ and a.e. $t\in\oi01$}, \label{eq: continuity eq on curve1}
			\\&\g'(t) = v(t,\g(t)) \text{ for $\cm$-a.e. $(\h,\g)$ and a.e. $t$ such that $\h(t)\neq0$},
			\\ \text{and }&\int_\HGF\int_0^1 (1+|v(t,\g(t))| + |g(t,\g(t))|)\,|\h(t)|\,\xdif t\xdif\abs{\cm}(\h,\g) <+\infty,\label{eq: continuity eq on curve3}
		\end{align}
		then $\int_\HGF\norm\h_\infty\xdif\abs\cm<+\infty$ and $\Theta(\cm)$ satisfies the continuity equation \eqref{eq: continuity equation}. 
		
		Conversely, given $\dm\in\MC(\X)$, if $\dm\geq0$ satisfies the continuity equation \eqref{eq: continuity equation} and 
		\begin{equation}
			\int_\X (1+|v(t,\x)|^2+|g(t,\x)|^2)\xdif\dm(t,\x) <+\infty,
		\end{equation}
		then $\dm=\Theta(\cm)$ for some $\cm\in\MC^+(\HGF)$ such that \eqref{eq: continuity eq on curve1}-\eqref{eq: continuity eq on curve3} hold and $\int_\HGF\norm\h_\infty\xdif\cm<+\infty$.
	\end{enumerate}
\end{thrm}
\edit{}{These results are mainly proved in \cite{Ambrosio2008,BrediesUnbalanced}, we extend to signed measures $\cm$ in the appendix (Theorem~\Rref{app: Theta properties}).} The results achieve two key relations: characterising when the mapping from $\MC(\X)$  to $\MC(\HGF)$ is well-defined, and when $\Theta(\cm)\in \CwX$. Weak continuity is of practical importance in applications. Without it, for example in the application of microscopy, we could not consider one frame of video to correspond to a single instance in time. Unfortunately we have seen that not all $\cm\in\MC(\HGF)$ satisfy this smoothness requirement. However, in the next lemma we will confirm that, if $\HG\subset \HGF$ is sufficiently ``small'', then $\Theta(\cm)\in\CwX$ for all $\cm\in\MC(\HG)$ due to the implicit assumption that $\norm\cm<+\infty$. A related property is the continuity of the operator $\Theta$ which we also confirm.
\begin{lmm}\label{lmma: Gamma_infty vs Gamma_1}
	For each $p\in\ci{1}{+\infty}$ define the set
	\begin{equation}
		\HG_p \eqdef \enscond{\hg=(\h,\g)\in\HGF}{\norm{\h}_p\leq 1},
	\end{equation}
	then
	\begin{equation}
		\enscond{\Theta(\cm)}{\cm\in\MC(\HGF),\ \int_{\HGF}\norm\h_p\xdif\abs\cm <+\infty} = \enscond{\Theta(\hat\cm)}{\hat\cm\in\MC(\HG_p)}
	\end{equation}
	and $\Theta\colon\MC(\HG_p)\to\MC(\X)$ is narrowly continuous. 
		
	Furthermore, if $p=+\infty$, then $\forall t\in\ci01$, $\cmt[\colon]{t}\MC(\HG_\infty)\to\MC(\bar\Omega)$ is also narrowly continuous. 
	
	In particular, sequentially we have that, for any sequence $\cm^n\xrightharpoonup{*}\cm$ narrowly in $\MC(\HG_p)$:
	\begin{align}
		\text{for all }p\in\ci{1}{+\infty}, &\qquad \Theta(\cm^n)\stackrel{*}{\rightharpoonup}\Theta(\cm) \text{ narrowly in }\MC(\X),
		\\\text{if }p=+\infty,\ \forall t\in\ci01,&\qquad \cmt[\cm^n]{t}\stackrel{*}{\rightharpoonup}\cmt{t} \text{ narrowly in }\MC(\bar\Omega).
	\end{align}
\end{lmm}
\edit{}{The proof of this lemma is found in Lemma~\Rref{app: Gamma_infty vs Gamma_1}, relying heavily on the definition of continuity given by \cite{Yosida1980}}. 
\begin{rmrk}
	Note that the definition of $\HG_p$ is consistent with the relation $\sim$, so $(\sfrac{\HG_p}{\sim},\metric)$ is also a metric space. In the rest of the paper, we require that $\HG$ is a closed subset of $\sfrac{\HG_\infty}{\sim}$ in order for it to be a complete separable metric space with $\Theta\colon\MC(\HG)\to\CwX$.
\end{rmrk}

Summarising the results of this section, in the remainder of this work we want to use a domain $D\subset\MC^+(\HGF)$ such that $\Theta|_D$ has nice properties with respect to the space $\CwX$. The combination of Theorem~\Rref{thm: Theta properties} and Lemma~\Rref{lmma: Gamma_infty vs Gamma_1} show that it is sufficient to consider either $D\subset\enscond{\cm\in\MC^+(\HGF)}{\int_{\HGF}\norm\h_\infty\xdif\cm<\infty}$ or simply $D\subset\MC^+(\HG_\infty)$. Analytically we will always consider $D\subset\MC^+(\HG_\infty)$ as it is more concise, although numerically either convention is equivalent. The generality of allowing any closed subset $\HG\subset\HG_\infty$ allows us to treat different applications with the same analysis, for example: 
\begin{description}
	\item[$\bf{\HG\subset\enscond{(\h,\g)\in\HG_\infty}{\h\equiv1}}$] This enforces balanced transport (\eg\ the Benamou-Brenier example \cite{BrediesBalanced}). If $\cm\in\MC^+(\HG)$, then $\Theta(\cm)\geq 0$ and mass is preserved on paths (\eg\ $t\mapsto\int_{\Omega}\xdif\Theta(\cm)_t$ is constant).
	\item[$\bf{\HG\subset\enscond{(\h,\g)\in\HG_\infty}{\h\in\Co(\ci01;\ci01)}}$] This allows unbalanced transport of non-negative mass (\eg\ the \\Wasserstein-Fisher-Rao example \cite{BrediesUnbalanced}). We still have $\Theta(\cm)\geq 0$, but $t\mapsto\int_{\Omega}\xdif\Theta(\cm)_t$ is not (necessarily) constant. In particular, mass can be created or destroyed (continuously) at any time.
	\item[$\bf{\HG\subset\HG_\infty}$] In the general case $\Theta(\cm)$ is a general signed measure, even when $\cm\geq 0$. In words, $\cm$ can give positive weight to curves with negative mass. The only constraint is that $\Theta(\cm)\in \CwX$ is continuous in time.
\end{description}

%% file: generalisation.tex

\section{Core variational problem}\label{sec: generalisation}
In this work we focus on inverse problems with dynamical but discrete-time structure. In particular, there exist observation times $t_j\in\ci01$, $j=0,\ldots,T$ and narrowly continuous linear operators $A_j\colon \MC(\bar\Omega)\to\xR^\m$. The operators $A_j$ are described by $a_i^j\in\Co(\bar\Omega)$ such that
\begin{equation}\label{eq: bounded linear kernels}
	\forall\dm \in\MC(\bar\Omega),\ i=1,\ldots,\m,\ j=0,\ldots,T, \qquad (A_j\dm)_i \eqdef \int_{\bar\Omega} a_i^j(\x)\xdif\dm(\x) .
\end{equation}
As stated at the end of Section~\Rref{sec: preliminary}, we work with a closed set of curves 
\begin{equation}
	\HG\subset\HG_\infty\eqdef\enscond{\hg=(\h,\g)}{\h\in\Co(\ci01;\ci{-1}{1}),\ \g\colon\ci01\to\bar\Omega,\ \g|_{\{\h\neq0\}} \text{ is continuous}},
\end{equation}
	so that $\forall t\in\ci01$, the map $\cmt[]{t}\colon\MC(\HG)\to\MC(\bar\Omega)$ is narrowly continuous. We therefore choose a data fidelity $\op{F}\colon\MC(\HG)\to\hc{0}{+\infty}$ of the form
\begin{equation}\label{eq: def of F}
	\text{for some convex }\op{F}_j\in \xCtwo(\xR^\m;\hc{0}{+\infty}),\qquad \op{F}(\cm) \eqdef \sum_{j=0}^T \op{F}_j(A_j \left[\cmt{t_j}\right]).
\end{equation}
For lower semi-continuous $w,\varphi\colon\HG\to\ci0{+\infty}$ define $\reg\colon\MC^+(\HG)\to\ho{-\infty}{+\infty}$, $D\subset\MC^+(\HG)$ by 
\begin{equation}
	\forall\cm\in D,\quad \reg(\cm) \eqdef \int_{\HG}w(\hg)\xdif\cm(\hg)
	 \qtext{where} D \eqdef \enscond{ \cm\in\MC(\HG) \quad}{\quad \cm\geq0,\quad \int_{\HG}\varphi(\hg)\xdif\cm \leq 1} \label{eq: D def}.
\end{equation}
We consider minimising the energy $\op{E}\colon D\to\ho{-\infty}{+\infty}$ defined by
\begin{equation}\label{eq: def of E}
	\forall\cm\in D,\qquad \op{E}(\cm) \eqdef \op{F}(\cm) + \reg(\cm).
\end{equation}
The motivation behind this in an Inverse Problems setting is that $\op{F}$ represents a smooth data fidelity with linear observations recorded at times $t_j$, and the combination of $\reg$ and $D$ represent a regularisation of the problem. The choice of $\varphi$ (hence of $D$) is often made to ensure the well-posedness of the model, but we are mostly interested in cases where the constraint  $\int_{\HG}\varphi(\hg)\xdif\cm \leq 1$ is not active, so that the choice of $\varphi$ has no impact on the set of minimisers.

\subsection{Existence of sparse minimisers}
Any choice of energy in this framework leads to a sparse reconstruction in the space of curves.
\begin{thrm}\label{thm: structure of E}
	If $A_j$ are given by \eqref{eq: bounded linear kernels} and $\varphi,w\colon\HG\to\ci{0}{+\infty}$ are lower semi-continuous, then $\op{F}$ and $\op{E}$ are lower semi-continuous. Furthermore, recall $\op{F}$ is bounded below. If $w$ or $\varphi$ have compact sub-levelsets, and $\inf_{\hg\in\HG}\varphi(\hg)>0$, then $\op{E}|_D$ has compact sub-levelsets. There exists a choice of minimiser $\cm^*\in \argmin_{\cm\in D} \op{E}(\cm)$ such that
	\begin{equation}\label{eq: sparse structure}
		\text{for some } a_i\geq0,\ \hg^i\in \HG,\qquad \cm^* = \sum_{i=1}^s a_i \delta_{\hg^i}
	\end{equation}
	for some $s\leq\m(T+1)+1$. If in addition $\int_\HG\varphi\xdif\cm^*<1$, then $s\leq \m(T+1)$.
\end{thrm}
The proof is given in Appendix~\Rref{app: structure of E}. The smaller value of $s$ will often be valid in practice, for example any choice $\varphi(\hg) \leq \frac{w(\hg)}{\op{E}(0)+1}$ is always sufficient.

\begin{rmrk}\label{rmk: BB energy on curves}
	The Benamou-Brenier example from Section~\Rref{sec:BB example intro} can be formulated in this setting with
	\begin{equation}\label{eq: BB def}
		\HG= \enscond{(\h,\g)\in\HGF}{ \h \equiv 1, \g\in\ACtwo(\ci01;\bar\Omega)},\quad \op{F}_j(A_j\dm) = \frac12\norm{A_j\dm-\data_j}_2^2, \qtext{and} w(\h,\g) = \int_0^1 \alpha + \frac\beta2|\g'(t)|^2 \xdif t.
	\end{equation}
	More details are given in Section~\Rref{sec:WFR example} where the Benamou-Brenier example is the limiting case $\delta\to+\infty$.  The choice of $\varphi$ in \cite{BrediesNumerics} was equivalent to $\varphi(\hg) = \frac{w(\hg)}{\op{E}(0)}$, whereas we suggest the default of $\varphi(\hg) = \frac{\alpha}{\op{E}(0)}$ which is easier to analyse. Both functions are strictly positive but sufficiently small to ensure $\int_\HG\varphi\xdif\cm^*<1$.
\end{rmrk}

\subsection{Discrete-time formulation}\label{sec: discrete time formulation}
Recall that $\HG\subset L^\infty(\ci01;\ci{-1}{1}\times\bar\Omega)$ is a space of continuous-time curves, denote the discrete-time space
\begin{equation}
	\tilde\HG \eqdef \enscond{(\hg(t_0),\ldots,\hg(t_T))}{\hg\in \HG}\subset (\ci{-1}{1}\times\bar\Omega)^{T+1}.
\end{equation}
Until now we have formulated $\op{E}$ as a function of measures $\cm\in\MC(\HG)$, but in this subsection we show that it can be thought of equivalently as a function $\tilde{\op{E}}$ of $\tilde\cm\in \MC(\tilde\HG)$. 
For the remainder of this section, without loss of generality, we assume there is a minimiser
\begin{equation}
	\cm^*\in \argmin\enscond{\op{F}(\cm) + \int_\HG w\xdif \cm}{\cm\in\MC^+(\HG)}
\end{equation}
which is also a minimiser of the energy considered in \eqref{eq: def of E}, \ie\ $\int\varphi\xdif \cm^*<1$. If that is not the case, one can use a Lagrange multiplier to form a modified energy with $w\gets w+\lambda \varphi$ for some $\lambda\geq0$. Our key observation is that $\cm^*$ must be supported on ``geodesics'' of $w$ (or $w+\lambda\varphi$) which interpolate the discrete-time curves.
\begin{lmm}\label{lmma: geodesic support}
	Suppose $w\colon\HG\to\xR\cup\{\infty\}$ is lower semi-continuous with compact sub-levelsets. Then, for all minimisers $\cm^*$ of \eqref{eq: lagrange mutiplier energy}, $\hg\in G(\hg(t_0),\ldots,\hg(t_T))$ for a.e. $\hg\in \op{supp}(\cm^*)$ where $G\colon\tilde\HG\rightrightarrows\HG$ is given by
	\begin{equation}
		G(\tilde\hg) \eqdef \argmin_{\hg\in\HG}\enscond{w(\hg)}{\forall j=0,\ldots,T,\ \hg(t_j)=\tilde\hg_j,\ w(\hg)<+\infty}.
	\end{equation}
\end{lmm}
The first part of the proof uses a measurable choice theorem from \cite{Brown1973}.
\begin{lmm}\label{lmma: measurable choice}
	There exists a measurable function $g\colon \HG\to\HG$ such that $g(\hg)\in G(\hg(t_0),\ldots,\hg(t_T))$ for all $\hg\in\HG$.
\end{lmm}
\begin{proof}[Proof of Lemma~\Rref{lmma: measurable choice}.]
	Let the set $U\subset \HG\times \HG$ be given by 
	\begin{equation}
		U\eqdef \enscond{(\hg^1,\hg^2)}{\hg^1\in\HG,\ \hg^2\in G(\hg^1(t_0),\ldots,\hg^1(t_T))}.
	\end{equation}
	Note that both $U$ and $G(\tilde\hg)$ are closed due to the lower semi-continuity of $w$. Finally, addressing the requirements of \cite[Thm. 1]{Brown1973}, we know that $(\sfrac\HG\sim,d_\HG)$ is a complete separable metric space, $U$ is a (closed) Borel set, and $G(\tilde\hg)$ is compact for each $\tilde\hg\in\tilde\HG$ (as it is either empty or a closed subset of a sub-levelset of $w$). We conclude the existence of $g$ as described in the claim.
\end{proof}
We can now return to the main result of this subsection.
\begin{proof}[Proof of Lemma~\Rref{lmma: geodesic support}.]
	The proof will show the contradiction that $\op{E}(g_\sharp\cm^*)<\op{E}(\cm^*)$ if $\cm^*$ is not strictly supported on the image of $G$. First we confirm that $\op{F}(g_\sharp\cm^*) = \op{F}(\cm^*)$. By \eqref{eq: bounded linear kernels}, for all $i=1,\ldots,m$, $j=0,\ldots,T$
	\begin{align}
		(A_j\cmt[(g_\sharp)\cm^*]{t_j}) = \int_{\bar\Omega}a_i^j(x)\xdif \cmt[(g_\sharp)\cm^*]{t_j}
		&= \int_\HG \h(t_j) a_i^j(\g(t_j)) \xdif\cm^*(\hg) \qtext{where} (\h,\g) = g(\hg)
		\\ &= \int_\HG \h(t_j) a_i^j(\g(t_j)) \xdif\cm^*(\hg) \qtext{where} (\h,\g) = \hg
		\\ &= (A_j\cmt[\cm^*]{t_j}).
	\end{align}
	Therefore each $\op{F}_j(A_j\cmt[(g_\sharp)\cm^*]{t_j}) = \op{F}_j(A_j\cmt[\cm^*]{t_j})$ as required. On the other hand, note $w(g(\hg))\leq w(\hg)$ always, suppose there exist $\epsilon,\delta>0$ such that $\cm^*(\enscond{\hg}{w(g(\hg))\leq w(\hg)-\delta})=\epsilon$. In which case
	\begin{equation}
		\reg(g_\sharp\cm^*) = \int_{w(g(\hg))\leq w(\hg)-\delta} w(g(\hg))\xdif\cm^* + \int_{w(g(\hg))> w(\hg)-\delta} w(g(\hg))\xdif\cm^* \leq \reg(\cm^*) - \delta\epsilon.
	\end{equation}
	Combining these equations, that would imply that $\op{E}(g_\sharp\cm^*)< \op{E}(\cm^*)$, contradicting the optimality of $\cm^*$. We conclude that $\cm^*(\enscond{\hg}{w(g(\hg))<w(\hg)})=0$ as required.
\end{proof}

This shows we can perform computations in the discrete-time space $\tilde\HG$ and later lift curves back to $\HG$ using the geodesics $G$. This holds both pointwise between $\tilde\HG\leftrightarrow\HG$ and with measures $\MC(\tilde\HG)\leftrightarrow\MC(\HG)$.
\begin{rmrk}\label{rmk: BB geodesics}
	For the Benamou-Brenier example, $w$ is given in Remark~\Rref{rmk: BB energy on curves}. If $\Omega$ is convex, then $G(\tilde\hg)=\{\hg\}$ where $\hg$ is the unique piecewise linear interpolant of the points $(t_j,\tilde\hg_j)$, as commented in \cite[Rem. 4.10]{BrediesNumerics}. Also,
	$$ \forall \tilde\hg=(\h_0,\g_0,\ldots,\h_T,\g_T)\in\tilde\HG,\qquad w(G(\tilde\hg)) = \alpha + \frac\beta2 \sum_{j=1}^T \frac{|\g_j-\g_{j-1}|^2}{t_j-t_{j-1}}. $$
	We therefore know that all minimisers $\cm^*$ are supported on the set
	\begin{equation}\label{eq: finite time BB}
		\HG = \enscond{(\h,\g) \in \Co(\ci01;\{1\}\times\bar\Omega)}{\g|_{\oi{t_{j-1}}{t_j}} \text{ is linear for each }j=1,\ldots,T}.
	\end{equation}
	Restricting to this domain of curves is much more computationally convenient without losing analytical accuracy.
\end{rmrk}

%% file: fwcvter.tex
\section{Frank-Wolfe convergence}\label{sec: FW convergence}
While Theorem~\Rref{thm: structure of E} guarantees an analytical structure of minimisers, we must now choose a reconstruction algorithm which is capable of taking advantage of this structure. As previously stated, we will use a variant of the Frank-Wolfe algorithm to take advantage of the sparse structure of reconstructions.

\subsection{The Frank-Wolfe algorithm}
\newcommand{\varx}{x}
\newcommand{\varX}{X}
\newcommand{\xmin}{x^*}
\newcommand{\fnrj}{f}
\newcommand{\cvx}{D}
\newcommand{\Xspace}{\mathcal{X}}
\newcommand{\exts}{s}
\newcommand{\extS}{S}
\newcommand{\dotp}[2]{\left\langle#1,#2\right\rangle}

In order to derive a variant of the Frank-Wolfe algorithm with inexact or stochastic steps, we first need to highlight a few properties which, to the best of our knowledge, have not been stated in the literature. We refer the reader to \cite{Jaggi2013} for a thorough introduction to the Frank-Wolfe algorithm. 
For the sake of generality and future reference, we work in an abstract setting, assuming that we want to minimize a convex function $\fnrj\colon \cvx\rightarrow \xR$, where $\cvx$ is a nonempty convex set of a Hausdorff locally convex vector space $\Xspace$. 

Standard results in convex analysis \cite[Chap. 7]{aliprantis_infinite_2006} ensure that for all $x \in \cvx$ and all $h \in \left(\cvx-\varx\right)$, the directional derivative
\begin{align*}
	\fnrj'(\varx;h)\eqdef	\lim_{t \to 0^+} \frac{\fnrj(\varx+th)-\fnrj(x)}{t}
\end{align*}
exists in $\xR\cup \{-\infty\}$ and is a convex function of $h$. The main steps of the Frank-Wolfe algorithm are given in Algorithm~\Rref{alg: abstract FW}. The standard choice \cite{frank-fw1956,demyanov-1970approximate} uses the \emph{Linear Minimisation Oracle} defined as
\begin{align}
	\exts^{n+1} =\LO(\varx^{n}) \in \argmin_{\exts \in \cvx} \left(\fnrj(\varx^{n})+\fnrj'(\varx^{n};\exts-\varx^{n})\right)= \argmin_{\exts \in \cvx} \fnrj'(\varx^{n};\exts-\varx^{n}).\label{eq: alg abstract lmo}
\end{align}
In the large majority of cases the existence of the $\LO$ is implied by $f$ being G\^ateaux-differentiable and $\cvx$ compact. 
The analysis in \cite{Jaggi2013} stresses that Algorithm~\Rref{alg: abstract FW} with the linear minimisation oracle \eqref{eq: alg abstract lmo} yields a \emph{minimising sequence}, \ie{} $\lim_{n \to +\infty} \fnrj(x^{n})= \inf_\cvx \fnrj$, provided the curvature of $\fnrj$ is finite,
\begin{equation}
	C_\fnrj\eqdef \sup_{\substack{\varx \in \cvx,\tilde\varx\in D\\\lambda\in\oi01}} \frac{\fnrj(\varx+\lambda(\tilde\varx-\varx)) - \fnrj(\varx)-\lambda \fnrj'(\varx;\tilde\varx-\varx) }{\lambda^2} <+\infty.\label{eq:alg curvature}
\end{equation}
In situations where the $\LO$ is not easy to compute, we can instead choose any $\exts^{n+1}\in\cvx$ and measure its suitability using the primal-dual gap
\begin{equation}\label{eq: def gap}
	\forall\varx,\exts\in\cvx,\qquad \gap(\varx;\exts) \eqdef \fnrj'(\varx;\exts-\varx) - \inf_{\tilde\exts\in\cvx} \fnrj'(\varx;\tilde\exts-\varx).
\end{equation}
In general $\gap\geq 0$ and $\gap(\varx;\LO(\varx))=0$, \edit{}{suppose $\gap(\varx^n;\exts^{n+1}) \leq \epsilon_n$ for some controlled error $\epsilon_n\geq0$. It has previously been shown in \cite[Thm. 1]{Jaggi2013} that $\varx^n$ is a minimising sequence whenever $\epsilon_n = O(1/n)$. 

In our context, the guarantee $\epsilon_n = O(1/n)$ would still be prohibitive for large $n$.} The linear minimisation oracle consists of finding a curve $\hg\in\HG$ which minimises a certain energy, see Section~\Rref{sec: fw back ip} below. \edit{}{Instead, we propose to use random discretisations of the domain so that implicitly $\liminf_{n\to+\infty}\epsilon_n=0$, without a guaranteed rate. In comparison to the result of \cite{Jaggi2013}, the relaxed assumption on $\epsilon_n$ is accounted for by the linesearch in Line~\Rref{alg: abstract linesearch} which implicitly selects large/small steps when $\epsilon_n$ is correspondingly large/small}, \ie\ when we have a lucky/unlucky draw. Our only assumptions on $\fnrj$ are that it is G\^ateaux differentiable with bounded curvature ($C_\fnrj < +\infty$) and $\gap(\varx;\exts) <+\infty$ for all $\varx,\exts\in\cvx$. Then, in the case that $\liminf_{n\to+\infty}\epsilon_n=0$ (possibly almost surely), we show that $\varx^n$ is a minimising sequence (almost surely).

\edit{}{One final aspect of our algorithm is the freedom which is granted in Line~\Rref{alg: abstract improve}: one may choose any point which has lower energy than the one provided by the linesearch. This is now a standard addition to the Frank-Wolfe algorithm to enable much faster practical convergence, see for instance \cite{brediesInverseProblemsSpaces2013,boydAlternatingDescentConditional2017,denoyelleSlidingFrankWolfeAlgorithm2019}. We exploit this in Section~\Rref{sec: results}.}

\begin{algorithm}
	\begin{algorithmic}[1]
		\State Choose $\varx^0\in \cvx$, $n\gets0$, 
		\Repeat
		\State Choose $\exts^{n+1} \in \cvx$ \Comment{oracle} \label{alg: abstract step}
		\State $\lambda_n\gets \argmin_{\lambda\in\ci01} \fnrj((1-\lambda)\varx^n + \lambda \exts^{n+1})$ \Comment{exact linesearch}\label{alg: abstract linesearch}
		\State Choose $\varx^{n+1}$ such that $\fnrj(\varx^{n+1}) \leq \fnrj((1-\lambda_n)\varx^n + \lambda_n \exts^{n+1})$ \Comment{improvement over the linesearch}\label{alg: abstract improve}
		\State $n\gets n+1$
		\Until{converged}
	\end{algorithmic}
	\caption{Abstract Frank-Wolfe algorithm}\label{alg: abstract FW}
\end{algorithm}

To ease the analysis of the randomised algorithm, we first analyse the deterministic version.

\begin{proposition}\label{prop: inexact FW}
	Let $\fnrj$, $\cvx$ be as above, and assume that $C_{\fnrj}<+\infty$ (see \eqref{eq:alg curvature}). If 
	\begin{align}
		\liminf_{n\to+\infty}\gap(\varx^n;\exts^{n+1})=0, \qtext{\ie} \liminf_{n\to +\infty} \left(\fnrj'(x^n;\exts^{n+1}-\varx^n) -\min_\cvx \fnrj'(\varx^n;\cdot-\varx^n) \right) = 0,\label{eq: alg liminf}
	\end{align}
	then Algorithm~\Rref{alg: abstract FW} yields a minimizing sequence.
\end{proposition}

\begin{proof}	
	First observe from Lines~\Rref{alg: abstract linesearch} and \Rref{alg: abstract improve}, and by definition of $C_{\fnrj}$ that
	\begin{align}
		\forall \lambda \in \ci{0}{1},\quad &\fnrj(\varx^{n+1})\leq \fnrj((1-\lambda)\varx^{n}+\lambda \exts^{n+1})\leq \fnrj(\varx^{n})+\lambda \fnrj'(\varx^{n};\exts^{n+1}-\varx^{n}) + C_\fnrj\lambda^2\\
		\mbox{hence}\quad &\fnrj(\varx^{n+1})-\fnrj(\varx^n) - \lambda \min_\cvx \left(\fnrj'(\varx^n;\cdot-\varx^n)\right) \leq \lambda \gap(\varx^n;\exts^{n+1}) + C_{\fnrj}\lambda^2. \label{eq: FW pre-liminf}
	\end{align}
	Note from the $\lambda=0$ case we have $\fnrj(\varx^n)\geq \fnrj(\varx^{n+1})$ for each $n$, therefore it is clear that $(\fnrj(\varx^n))_{n\in\xN}$ converges to some limit $\ell\geq \inf_{\varx\in\cvx}\fnrj(\varx) \geq -\infty$. 
	We are required to show that $\ell=\inf_{\varx\in\cvx}\fnrj(\varx)$. 
	
	The case $\ell=-\infty$ is trivial, otherwise we also have $\lim_{n\to+\infty} \fnrj(\varx^{n+1})-\fnrj(\varx^{n}) = 0$. Now, taking the $\liminf$ on both sides of \eqref{eq: FW pre-liminf} gives
	\begin{equation}
		\liminf_{n\to+\infty}-\lambda \min_\cvx \left(\fnrj'(\varx^{n};\cdot-\varx^{n})\right) \leq C_{\fnrj}\lambda^2.
	\end{equation}
	Dividing by $\lambda\to0^+$, we obtain $\limsup_{n\to+\infty}\min_\cvx \left(\fnrj'(\varx^{n};\cdot-\varx^{n})\right) \geq 0$. On the other hand, by convexity,
	\begin{align}
	\forall \varx \in \cvx,\quad \fnrj(\varx) \geq \limsup_{n\to+\infty} \left[\fnrj(\varx^{n}) + \fnrj'(\varx^{n};\varx-\varx^{n})\right] 
		\geq \ell + \limsup_{n\to+\infty} \min_\cvx \left(\fnrj'(\varx^{n};\cdot-\varx^{n})\right)
		\geq \ell.
	\end{align}
	Since $x \in \cvx$ is arbitrary, we deduce that $\ell=\lim_{n \to \infty} \fnrj(\varx^n)= \inf_\cvx \fnrj$, as required.
\end{proof}

\subsection{Stochastic variant}
Now, we may study the behaviour of the algorithm in a stochastic framework. We build a random process $(\varX^{n},\extS^n)_{n\in \xN^*}$ in $\cvx\times \cvx$, by considering a random initialisation $\varX^0$ in $\cvx$, and applying Algorithm~\Rref{alg: abstract FW}  by picking a random point $\extS^{n+1}$ in $\cvx$ at Line~\Rref{alg: abstract step}. Note that, at each step, $\extS^{n+1}$ may depend on $\{\varX^k\}_{k=0}^{n}$ and $\{\extS^{k}\}_{k=1}^{n}$.

Typically, as in Section~\Rref{sec: fw back ip}, we consider a setting where solving \eqref{eq: alg abstract lmo} exactly is too costly, and where one draws a random grid on which to perform the optimisation. The variable $\extS^{n+1}$ is then a minimizer of $\fnrj'(\varX^{n};\cdot-\varX^n)$ among a finite, small, subset of $\cvx$.

Let $(\Filtration^n)_{n \in \xN}$ be the filtration generated by that random process, \ie{} $\Filtration^n$ is the $\sigma$-algebra generated by $\{\varX^0\}\cup \{\varX^k\}_{1\leq k\leq N}\cup \{\extS^k\}_{1\leq k\leq N}$.

\begin{proposition}\label{prop: liminf}
	Let $(\varX^{n})_{n\in \xN}$ and  $(\extS^{n})_{n \in \xN^*}$ as described above, and let $(\Filtration^n)_{n \in \xN}$  be the filtration they generate. 

	If for all $\epsilon>0$, all $n \in \xN$, there is some deterministic $p_n(\epsilon)>0$ such that 
\begin{align*}
	\P( \gap(\varX^n;\exts^{n+1}) < \epsilon | \Filtration^n)\geq p_n(\epsilon)\quad  \mbox{almost surely,}
\end{align*}
and $\sum_{n=1}^{\infty} p_n(\epsilon)=+\infty$,
then $(\fnrj(\varX^n))_{n\in \xN}$ is a minimizing sequence almost surely.

\end{proposition}
The above proposition is a consequence of the following lemma, setting $G^{n+1}=\gap(\varX^n;\exts^{n+1})$.

\begin{lmm}\label{lmma: a.s. FW convergence}
	Let $(\Filtration^{n})_{n\in \xN}$ be a filtration, and $(G^{n})_{n\in \xN}$ a family of random variables such that for each $n,m \in \xN$, $n\leq m$, $G^{n}$ is $\Filtration^{m}$-measurable. If for all $\epsilon>0$, $n \in \xN$, there is some deterministic $p_n(\epsilon)\geq 0$ such that 
	\begin{align}
		\P(G^{n+1} < \epsilon | \Filtration^n)\geq p_n(\epsilon)\quad  \mbox{almost surely}
	\end{align}
	and $\sum_{n=0}^{\infty} p_n(\epsilon)=+\infty$, then $\liminf_{n \to \infty} G^n \leq 0$ almost surely.	
\end{lmm}

\begin{proof}
	Fix $\varepsilon>0$. For all $N, M \in \xN$ with $M\geq N$, 
	\begin{align}
		\P\left(\bigcap_{n\geq N}\{G^n\geq \varepsilon\}\right)\leq \P \left(\bigcap_{n=N}^{M+1}\{G^n\geq \varepsilon\}\right)&= \E \left[ \1_{\{G^{M+1}> \varepsilon\}} \left(\prod_{n=N}^{M}\1_{\{G^n\geq \varepsilon\}}\right)  \right]\\
		&\leq \E \left[ \E\left(\1_{\{G^{M+1}\geq \varepsilon\}} | \Filtration^M\right) \prod_{n=N}^{M}\1_{\{G^n\geq \varepsilon\}}\right]\\
		&\leq (1-p_{M+1}(\varepsilon))\P\left(\bigcap_{n\geq N}\{G^n\geq \varepsilon\}\right)\\
		&\leq \prod_{n=N}^{M+1}(1-p_n(\varepsilon))\\
		&\leq \exp\left(-\sum_{n=N}^{M+1} p_n(\varepsilon) \right).
	\end{align}
	
	Letting $M\to +\infty$, we get $\P\left(\bigcap_{n\geq N}\{G^n\geq \varepsilon\}\right)=0$, that is $\P\left(\bigcup_{n\geq N}\{G^n< \varepsilon\}\right)=1$. As a result,
	
	\begin{align*}
		\P\left(\liminf_{n \to \infty} G^n\leq 0\right)= \P\left(\bigcap_{k \in \xN^*}\bigcap_{N \in \xN}\bigcup_{n\geq N} \left\{G^n< \frac{1}{k}\right\} \right)= \lim_{k \to +\infty} \lim_{N \to +\infty} \P\left(\bigcup_{n\geq N} \left\{G^n< \frac{1}{k}\right\} \right)=1.
	\end{align*}
\end{proof}
To summarise, the main convergence requirement of Proposition~\Rref{prop: inexact FW} is to guarantee $\liminf_{n\to+\infty} G^n \leq 0$. Our solution to this is to use random discretisations so that the stochastic variant of Algorithm~\Rref{alg: abstract FW} still converges asymptotically (almost surely), but the complexity of each individual iteration remains low. \edit{}{A similar idea in the setting of stochastic Frank-Wolfe was pursued in \cite{Silveti2020} where their proof relies on what is called Assumption P.8, in our notation this requires the sum $\sum_{n=0}^\infty\lambda_n G^n<+\infty$ to be finite almost surely, where $\lambda_n$ is chosen deterministically. To apply this algorithm in our setting we would therefore need to bound the magnitude of $G^n$ relative to the \apriori\ choice of $\lambda_n$. With Proposition~\Rref{prop: liminf}, we overcome this limitation with the linesearch for $\lambda$}, so the only remaining requirement is to guarantee a uniform probability of achieving a good search direction.

\subsection{Back to the dynamic inverse problem}\label{sec: fw back ip}
We consider again the setting of Section~\Rref{sec: generalisation}. 
The function $\op{E}$ is convex on $D$, and we note that, for each $\cm \in D$, its directional derivative is given by 
\begin{align*}
	\forall \nu \in D,\quad	\op{E}'(\cm;\nu-\cm) = \int_{\HG}\left(\op{F}'(\cm)(\hg) + \reg'(\cm)(\hg)\right)\xdif (\nu-\cm)(\hg),
\end{align*}
with $\reg'(\cm)=\left[\hg\mapsto w(\hg)\right]$ and 
\begin{align*}
	\op{F}'(\cm)= \left[(\h,\g)\mapsto\sum_{j=0}^T \h(t_j)\eta_j(\g(t_j))\right]\in \Cb(\HG),\qtext{where} \eta_j \eqdef A_j^*\nabla \op{F}_j(A_j\cmt{t_j}).
\end{align*}
In particular, \edit{}{note that we can  write $\op{E}'(\cm;\ep-\cm) = \int_\HG \op{E}'(\cm)\xdif[\ep-\cm]$ with $\op{E}'(\cm)\colon\HG\to\xR$.} This structure enables us to prove that $\op{E}$ satisfies the major requirement of Section~\Rref{sec: FW convergence}, that $\inf_{\ep\in D} \op{E}'(\cm;\ep-\cm)>-\infty$ for all $\cm\in D$. To do this, we show that the infimum is always achieved by an extreme point of $D$.

\begin{lmm}\label{lmma: well defined linear oracle}
	Suppose $\varphi,w\colon\HG\to\ci{0}{+\infty}$ are lower semi-continuous. If $\inf_{\hg\in\HG}\varphi(\hg)>0$ and either $\varphi$ or $w$ have compact sub-levelsets, then $\op{E}'(\cm;\cdot-\cm)$ is lower semi-continuous and coercive on $D$. Moreover it has minimiser of the form
	\begin{equation}\label{eq: exact linear oracle}
		\LO(\cm) \in\{0\} \cup \enscond{\varphi(\hg^*)^{-1}\delta_{\hg^*}}{ \hg^* \in \argmin_{\hg\in\HG} \frac{\eta(\hg) + w(\hg)}{\varphi(\hg)}}.
	\end{equation}
	where $\eta(\h,\g) \eqdef \sum_{j=0}^T \h(t_j)\eta_j(\g(t_j))\in \Cb(\HG)$.
\end{lmm}
\begin{proof}
	Recall the definition of $\LO$, 
	\begin{equation}\label{eq: E tilde def}
		\LO(\cm) \in \argmin_{\tilde\cm\in D} \tilde{\op{E}}(\tilde\cm) \qtext{where} \tilde{\op{E}}(\tilde\cm) \eqdef \int_{\HG}\left(\eta +w\right)\xdif\tilde\cm.
	\end{equation}
	The properties of $D$ come from Lemma~\Rref{lmma: extreme point computation}, in particular $D$ is convex, closed and bounded, so $\inf_{\tilde\cm\in D}\int_\HG\eta\xdif\tilde\cm>-\infty$ and $\tilde{\op{E}}$ is lower semi-continuous (Lemma~\Rref{thm: Fatou's lemma}). It is therefore well-posed to consider minimisers of $\tilde{\op{E}}$. We show that there is a choice $\LO(\cm)=\cm^*\in\Ext D$.
	
	\begin{description}
		\item[Case $D$ is compact] If $\varphi$ has compact sub-levelsets, then $D$ is compact by Lemma~\Rref{lmma: extreme point computation}. Bauer's principle therefore states that there exists a point
		\begin{equation}
			\cm^*\in\Ext{D} \qtext{such that} \tilde{\op{E}}(\cm^*) = \inf_{\tilde\cm\in D}\tilde{\op{E}}(\tilde\cm).
		\end{equation}
		
		\item[Else] Otherwise, $w$ has compact sub-levelsets, so the sub-levelset 
		\begin{equation}
			U \eqdef \enscond{\tilde\cm\in D\quad}{\quad \tilde{\op{E}}(\tilde\cm) \leq 1}
		\end{equation}
		is compact by Theorem~\Rref{thm: compact sub-levelsets}, convex, and non-empty because $0\in D$, $\tilde{\op{E}}(0)=0$. Application of Bauer's principle now gives a point $\cm^*\in\Ext{U}$, to complete the proof we will confirm $\cm^*\in\Ext{D}$. 
		
		Suppose for contradiction that there exists $\cm_0,\cm_1\in D\setminus\{\cm^*\}$ with $\cm^*\in\oi{\cm_0}{\cm_1}$. Without loss of generality, we assume that $\cm^*=\frac{1}{2}\cm_0+\frac{1}{2}\cm_1$ and we define
		\begin{equation}
			\forall \lambda\in\ci01,\qquad \cm_\lambda\eqdef\lambda\cm_1 + (1-\lambda)\cm_0\in D,\qquad \cm_{\frac12} = \cm^*.
		\end{equation}
		Consider the function $f(\lambda) = \tilde{\op{E}}(\cm_\lambda)$, which is linear and $f(\frac12)=\tilde{\op{E}}(\cm^*) \leq \tilde{\op{E}}(0)=0$. As $f$ is continuous, there exists $\epsilon>0$ with $f(\frac12\pm\epsilon)\leq \frac12$. We conclude that $\cm_{\frac12\pm\epsilon}\in U$ and $\cm^*\in\oi{\cm_{\frac12-\epsilon}}{\cm_{\frac12+\epsilon}}$, contradicting the assumption that $\cm^*\in\Ext{U}$.
	\end{description}
	In both cases we see there is a choice $\LO(\cm)=\cm^*\in\Ext D$. Finally, by Lemma~\Rref{lmma: extreme point computation} we have 
	\begin{equation}
		\Ext D = \{0\}\cup\enscond{\varphi(\hg)^{-1}\delta_{\hg}}{\varphi(\hg)<+\infty},
	\end{equation} 
	so we deduce \eqref{eq: exact linear oracle} as required.
\end{proof}

This lemma is also useful for applying the stochastic Frank-Wolfe algorithm, Lemma~\Rref{lmma: a.s. FW convergence}. For each $n\in\xN$, let $z^n\sim Z$ be an independent random discretisation of the domain $\tilde\HG = (\ci{-1}1\times\bar\Omega)^{T+1}$ (following the discrete-time formulation in Section~\ref{sec: discrete time formulation}). Continuing with the notation $\tilde{\op{E}} = \op{E}'(\cm^{n-1})$ from \eqref{eq: E tilde def}, we will use the discrete minimiser denoted 
\begin{equation}
	\ep^n = \bar{\LO}(\cm^{n-1},z^n) \in \{0\}\cup \enscond{\varphi(\hg_d^*)^{-1}\delta_{\hg_d^*}}{\hg_d^*\in \argmin_{\hg\in z^n} \tilde{\op{E}}\left(\varphi(\hg)^{-1}\delta_\hg\right)}.
\end{equation}
The remaining requirements for applying both Frank-Wolfe variants simplify greatly in the Benamou-Brenier example where $\varphi>0$ is a constant, curves have constant mass, and $\Omega=\oi01^d$ (\ie\ $z^n\subset \ci01^{d(T+1)}$). An immediate consequence of this is uniform boundedness, as $\int_\HG\xdif\cm=\sum_{j=0}^T\norm{\cmt{t_j}}\leq \varphi^{-1}$ for all $\cm\in D$. Because $\op{E}$ is quadratic and $A_j$ bounded (see \eqref{eq: bounded linear kernels}, $(A_j\dm)_i = \int_{\bar\Omega}a_i^j\xdif\dm $), the curvature bound follows immediately:
\begin{equation}
	C_{\op{E}} = \sup_{\cm,\cm'\in D} \sum_{j=0}^T\norm{A_j\cmt[{[\cm-\cm']}]{t_j}}_2^2 
	\leq \left(\sup_{\cm,\cm'\in D} \sum_{j=0}^T\norm{A_j}\norm{\cmt[{[\cm-\cm']}]{t_j}}\right)^2
	\leq \varphi^{-2}\max_{i,j}\norm{a_i^j}_\infty < +\infty.
\end{equation}
Finally we must show that $\gap(\cm^{n-1};\ep^n)$ is uniformly small, independently of $\cm^{n-1}$. Ignoring the $\cm^*=0$ case, this quantity can be written
\begin{equation}
	\gap(\cm^{n-1};\ep^n) = \int_\HG \tilde{\op{E}}\xdif [\ep^n - \LO(\cm^{n-1})] \leq \varphi^{-1}\min_{\hg_d^*\in z^n}\max_{\hg^*\in \HG} [\tilde{\op{E}}(\hg^*_d) - \tilde{\op{E}}(\hg^*)]
\end{equation}
where, from \eqref{eq: BB def} and \eqref{eq: finite time BB}, for all $\g\in \bar\Omega^{T+1}$ we have
\begin{equation}\label{eq: BB specific LMO}
	\tilde{\op{E}}(1,\g) = \left[\sum_{i=1}^\m\sum_{j=0}^T \left(\int_\Omega a_i^j(x)\xdif\cmt[\cm^{n-1}]{t_j}(x) - b_j\right)^\top a_i^j(\g_j)\right] + \alpha + \frac\beta2\sum_{j=1}^T \frac{|\g_j-\g_{j-1}|^2}{t_j-t_{j-1}}
\end{equation}
for some $\alpha,\beta>0$. If moreover $A_j$ are Lipschitz, then $\tilde{\op{E}}$ is also uniformly Lipschitz independently of $\cm^{n-1}$. Combining this uniform Lipschitz property with, for example, uniformly sampled discrete grids $z^n$ guarantees the uniform bound for Lemma~\Rref{lmma: a.s. FW convergence}.

\subsection{Choice of constraint}\label{sec:FW discussion}
In this section we will briefly discuss the main difference between our formulation of the Benamou-Brenier example and that implemented in \cite{BrediesNumerics}. Although the construction in \cite{BrediesNumerics} looks very different, it can also be seen as Frank-Wolfe/Generalised Conditional Gradient approach to minimise the same function $\op{E}$ (this equivalence is confirmed in Section~\Rref{sec:WFR example}). The parallel result to Lemma~\Rref{eq: exact linear oracle} is \cite[Prop. 3.6]{BrediesNumerics} which shows that they are also incrementally adding new curves $\g$ to the support of the reconstruction $\cm$ each iteration.  Both implementations actually use a stochastic variant of Frank-Wolfe, but we will discuss the classical version for simplicity.

Recall from Remark~\Rref{rmk: BB energy on curves}, that the un-constrained energy is
\begin{align}
	\forall \cm \in \MC^+(\HG),\quad \op{E}(\cm) = \frac12\sum_{j=0}^T\norm{A_j\cmt{t_j}-\data_j}_2^2 
	+\int_\HG w\xdif\cm\label{bbexample-eq-objfw}
\end{align}
where $w(\g)= \alpha + \frac{\beta}{2}\int_0^1|\g'(t)|^2\xdif t$ for all $\g\in\HG$. At iteration $n$, let $\eta(\g) \eqdef \sum_{i=1}^\m\sum_{j=0}^T \left(A_j\cmt[\cm^n]{t_j}-\data_j\right)_ia_i^j(\g(t_j))$ denote the linearisation of the fidelity term where operators $A_j$ are represented by kernels $a_i^j\in\Cb(\bar\Omega)$. The standard Frank-Wolfe procedure consists in iteratively minimising the function 
\begin{equation}\label{eq: BB linear oracle}
	\min_{\cm \in D}\int_{\HG} \left[\eta(\g) + w(\g)\right]\xdif\cm.
\end{equation}
on a constraint set $D \subset \MC^+(\HG)$ to yield descent directions. 
Different sets $D\supset \argmin_{\cm\in\MC^+(\HG)}\op{E}(\cm)$ can change the minimisers of \eqref{eq: BB linear oracle} without necessarily changing the original problem \eqref{bbexample-eq-objfw}. A first choice is to consider
\begin{equation}
	D_1\eqdef \enscond{\cm \in \MC^+(\HG)}{\int_{\HG} w\xdif\cm\leq \op{E}(0)},
\end{equation}
which does contain the minimisers of \eqref{bbexample-eq-objfw} since they must satisfy $\int_\HG w\xdif \cm\leq \op{E}(\cm)\leq \op{E}(0)$.
By Lemma~\Rref{lmma: extreme point computation}, the set of extreme points of $D_1$ is 
$\{0\}\cup \enscond{ \frac{\op{E}(0)}{w(\g)}\delta_{\g}}{\g \in \HG}$, which yields the descent direction 
\begin{equation}\label{eq: BB Bredies LO}
	\g^{n+1} \in \argmin_{\g\in\HG} \frac{\eta(\g) + w(\g)}{w(\g)} = \argmin_{\g\in\HG} \frac{\eta(\g)}{w(\g)}.
\end{equation}
That is precisely the descent direction used in~\cite[(4.30)]{BrediesNumerics}.
On the other hand, we propose to use
\begin{equation}
	D_2\eqdef \enscond{\cm \in \MC^+(\HG)}{\int_{\HG} \alpha \xdif\cm\leq \op{E}(0)},
\end{equation}
which also contains the minimisers of \eqref{bbexample-eq-objfw}, since they must satisfy ${\int_\HG \alpha\xdif \cm \leq \int_\HG w\xdif \cm\leq \op{E}(\cm)\leq \op{E}(0)}$.
Applying Lemma~\Rref{lmma: extreme point computation}, we see that the set of extreme points of $D_2$ is 
$\{0\} \cup \enscond{\frac{\op{E}(0)}{\alpha}\delta_{\g}}{\g \in \HG}$, which yields 
\begin{equation}\label{eq: BB our LO}
	\g^{n+1} \in \argmin_{\g\in\HG} \frac{\eta(\g) + w(\g)}{1} = \argmin_{\g\in\HG} \eta(\g) + w(\g).
\end{equation}
We found \eqref{eq: BB our LO} more convenient than \eqref{eq: BB Bredies LO} as it is amenable to dynamic programming techniques (see Section~\Rref{sec: dynamical programming}). The whole minimisation algorithm is summarised in Algorithm~\Rref{alg: FW}.

\begin{rmrk}
	In the classical Frank-Wolfe algorithm $D$ must be compact. In measure spaces, compactness is equivalent to boundedness and tightness (Prokhorov's theorem). In both cases, boundedness comes from $\alpha>0$, and tightness comes from the fact that $w$ has compact sub-levelsets. The set $D_1$ is compact by Lemma~\Rref{thm: compact W sublevel}. On the other hand, $D_2$ is only bounded, but sub-levelsets of \eqref{eq: BB linear oracle} are still compact, as argued in Lemma~\Rref{lmma: well defined linear oracle}.
\end{rmrk}

\begin{algorithm}
	\begin{algorithmic}[1]
		\State Choose $ \cm^0=0$, $n\gets0$
		\Repeat
		\State Compute $\g^{n+1}$ according to \eqref{eq: BB Bredies LO} or \eqref{eq: BB our LO}\Comment{Linear oracle step} \label{alg: BB linear oracle line}
		\State Choose $\ep^{n+1} \propto \delta_{\g^{n+1}} $ such that $\ep^{n+1}\in\Ext{D} $ \Comment{Choice of scaling}
		\State $\cm^{n+1} \gets (1-\lambda_n)\cm^n + \lambda_n \ep^{n+1}$ \Comment{Some stepsize $\lambda_n\in\ci01$}
		\State $n\gets n+1$
		\Until{converged}
	\end{algorithmic}
	\caption{Frank-Wolfe algorithm for the Benamou-Brenier example}\label{alg: FW}
\end{algorithm}

%% file: generaldp.tex
\section{Dynamical programming method}\label{sec: dynamical programming}
Throughout Section~\Rref{sec: FW convergence} we have shown that we can find a minimising sequence to $\op{E}$ from Section~\Rref{sec: generalisation} by repeatedly evaluating a simplified linear oracle. In particular, we compute the minimiser from Lemma~\Rref{lmma: well defined linear oracle}, but over a discretised domain. In Section~\Rref{sec: discrete time formulation} we motivated discretising in time to $\tilde\HG=(\ci{-1}{1}\times\bar\Omega)^{T+1}$ without loss of precision, now we also discretise in space using the domain $\mesh \eqdef \prod_{j=0}^T\mesh_j$ for some discrete sets $\mesh_j \subset \ci{-1}{1}\times\bar\Omega$ (``grids'' in the mass-location space). Whilst $\mesh$ is much smaller than $\tilde\HG$, naive computation time for $\hg^*$ still scales exponentially in $T$. We propose to compute this discrete minimiser efficiently using dynamical programming, for which we require two final simplifications:
\begin{align}
	\varphi(\hg) &= \varphi_0 &&\text{ for some constant }\varphi_0>0, \text{ and}\tag{A1}\label{eq: const constraint}
	\\w(\hg) &= \sum_{j=1}^T\step_j(\hg(t_{j-1}),\hg(t_j)) &&\text{ for some }\step_j\colon (\ci{-1}{1}\times\bar\Omega)^2\to\xR. \tag{A2}\label{eq: step penalty}
\end{align}
\begin{rmrk}
	For the Benamou-Brenier example, these assumptions are satisfied by the choice 
	\begin{equation}
		\varphi_0 = \frac{\alpha}{\op{E}(0)},\qquad \step_j(\hg(t_{j-1}), \hg(t_{j})) = \alpha (t_j-t_{j-1}) + \frac\beta2\frac{|\g(t_j)-\g(t_{j-1})|^2}{t_j-t_{j-1}}
	\end{equation}
	for all $\hg= (\h,\g) \in \HG$. 
	We will see later in Section~\Rref{sec:WFR discrete-time energy} that these assumptions are also satisfied in the Wasserstein-Fisher-Rao example. In particular, the step function of the Benamou-Brenier penalty is Wasserstein optimal transport, and the step function of the dynamic Wasserstein-Fisher-Rao penalty is the static Wasserstein-Fisher-Rao penalty.
\end{rmrk}
Under assumptions \eqref{eq: const constraint} and \eqref{eq: step penalty}, the discretised optimisation problem \eqref{eq: exact linear oracle} can be greatly simplified to
\begin{align}
	\hg^* &\in \argmin_Y\enscond{\sum_{j=0}^T\eta_j(Y_j) + \sum_{j=1}^T\step_j(Y_{j-1},Y_j)}{Y\in \prod_{j=0}^T \mesh_j} 
	\\&\qtext{where} \forall j=0,\ldots,T,\qquad\eta_j(\h_j,\g_j) \eqdef \h_j [A_j^*\nabla \op{F}_j(A_j\cmt{t_j})](\g_j).
\end{align}
This minimisation problem can now be formulated as computing the minimal path on a weighted, directed, acyclic graph. The vertices are the points $(t_j,y)$ for $y\in \mesh_j$, the edge weights are given by $\eta_j$ and $\step_j$, and the time index provides an ordering and prevents cycles in the graph. Algorithms for computing minimal paths on directed acyclic graphs are well studied, the complexity bound is given below.
\begin{thrm}[{\cite[Sec. 24.2]{Cormen2009}}]\label{thm: linear oracle complexity}
	Suppose assumptions \eqref{eq: const constraint} and \eqref{eq: step penalty} hold and $|\mesh_j|\leq N$ for each $j$, then $\hg^*$ can be computed with complexity
	\begin{equation}\label{eq: search complexity}
		\text{total time} = O\left(N\sum_{j=0}^T\op{cost}(\eta_j) + N^2\sum_{j=1}^T\op{cost}(\step_j)\right)
	\end{equation}
	where $\op{cost}(\eta_j)$ is the cost of evaluating $\eta_j$ once \etc.
\end{thrm}
\begin{rmrk}
	In graph terminology, the two terms of \eqref{eq: search complexity} represent the number of vertices plus the number of edges respectively, asymptotically $O(NT)$ and $O(N^2T)$ respectively. If, for example, the paths have a maximum velocity $V$, then the number of edges per vertex is reduced from $N$ to $O(N(VT^{-1})^d)$. This results in a reduced total complexity of $O(NT + N^2T^{1-d}V^d)$.
\end{rmrk}

Finally we outline how the minimal path can be computed efficiently in our specific example by using dynamic programming. To do so, define the truncated energies and minimal paths:
\begin{align}
	\tilde{\op{E}}_J(Y) &\eqdef \sum_{j=0}^J \eta_j(Y_j) + \sum_{j=1}^J\step_j(Y_{j-1},Y_j),
	\\ Y^*[y,J] &\in \argmin_Y\enscond{\tilde{\op{E}}_{J}(Y)}{Y \in \prod_{j=0}^J\mesh_j,\ Y_J=y}
\end{align}
for each $J=0,\ldots,T$ and $y\in\mesh_J$. This generates $\hg^*$ through the computation
\begin{equation}
	\hg^* = Y^*[y^*,T], \qquad y^* \in \argmin_{y\in\mesh_T}\tilde{\op{E}}_T(Y^*[y,T]).
\end{equation}
Observe that for all $J= 1,\ldots,T$ and $y \in \mesh_J$,
\begin{align}
	\min_{\substack{Y \in \prod_{j=0}^J\mesh_j\\ Y_J=y}} \tilde{\op{E}}_J(Y) &= \min_{\substack{Y \in \prod_{j=0}^{J}\mesh_j}} \left\{\tilde{\op{E}}_{J-1}(Y) + \eta_{J}(y) + \step_J(Y_{J-1},y) \right\}\\
	&=\eta_{J}(y) + \min_{\substack{y' \in \mesh_{J-1}}}\Bigg\{\step_J(y',y) + \min_{\substack{Y \in \prod_{j=0}^{J-1}\mesh_j\\ Y_{J-1}=y'}} \tilde{\op{E}}_{J-1}(Y)\Bigg\}.
\end{align}
In particular, we can choose $Y^*[y, J]$ inductively to be
\begin{equation}\label{eq: inductive extreme point computation}
	Y^*[y,J]=
	\begin{pmatrix}
		Y^*[y',J-1] & y
	\end{pmatrix} \qtext{for} y'\in\argmin_{y' \in \mesh_{J-1}} 
	\left[\step_J(y',y) + \tilde{\op{E}}_{J-1}(Y^*[y',J-1])\right].
\end{equation}
For each $y$ and $J$ each of these steps requires $O(N)$ computation, confirming the global complexity of $O(N^2T)$.

%% file: WFRexample.tex

\section{The unbalanced Wasserstein-Fisher-Rao example}\label{sec:WFR example}
The two numerical examples considered in this work use the Benamou-Brenier and Wasserstein-Fisher-Rao (WFR) penalties. The former has already been discussed in previous remarks and is a limiting case of the WFR penalty so we will not discuss it in further detail here. 

The penalty considered in \cite{BrediesUnbalanced} is $\creg\colon\MC^+(\X)\to\xR\cup\{+\infty\}$ such that for all $\dm\in\MC^+(\X)$,
\begin{equation}\label{eq:WFR energy}
	\creg(\dm) \eqdef \inf_{\substack{v\in \xLtwo_\dm(\ci01\times\bar\Omega;\xR^d) \\ g\in \xLtwo_\dm(\ci01\times\bar\Omega)}}\Bigg\{\int_\X \left[\alpha + \frac\beta2|v|^2+\frac{\beta\delta^2}2g^2\right]\xdif\dm \ \st\ \partial_t\dm + \op{div}(v\dm)=g\dm\Bigg\}
\end{equation}
for some $\alpha,\beta,\delta>0$, and the continuity equation is satisfied in the sense of \eqref{eq: continuity equation}. This leads to the energy
\begin{equation}\label{wfrexample-eq-energy}
	\C{E}(\dm) = \frac12\sum_{j=0}^T\norm{A_j\dm_{t_j}-\data_j}_2^2 + \creg(\dm)
\end{equation}
where the properties of $A_j$ are as stated in \eqref{eq: bounded linear kernels}.
Much is already known about minimisers of this energy.
\begin{thrm}[{\cite[Thms. 4.2,6.4]{BrediesUnbalanced}}]\label{thm: WFR sparsity}
	Let $\alpha,\beta,\delta>0$ and
	\begin{equation}
		\HG \eqdef \enscond{(\h,\g)\colon\ci01\to\ci01\times\bar\Omega}{\sqrt\h\in\ACtwo(\ci01),\ \sqrt\h\g\in\ACtwo(\ci01;\xR^d)}.
	\end{equation}
	Then for all $\dm\in\MC^+(\X)$ with $\creg(\dm)<+\infty$, there exists $\cm\in \MC^+(\HG)$, $v\in \xLtwo_\dm(\X;\xR^d)$, $g\in \xLtwo_\dm(\X)$ such that $\dm=\Theta(\cm)$, 
	\begin{equation}
		\creg(\dm) = \int_\X \left[\alpha + \frac\beta2|v|^2 + \frac{\beta\delta^2}2g^2\right]\xdif\dm,\qquad \partial_t\dm + \op{div}(v\dm)=g\dm,
	\end{equation} 
	and
	\begin{equation}\label{eq: WFR v = x'}
		\begin{aligned}
			\g'(t) &= v(t,\g(t))\qtext{for} \text{a.e. } t\in\{\h>0\}\text{ and } \cm\text{-a.e. }(\h,\g),
			\\\h'(t) &= g(t,\g(t))\h(t)\qtext{for} \text{a.e. } t\in\ci01 \text{ and } \cm\text{-a.e. }(\h,\g).
		\end{aligned}
	\end{equation}	
	Moreover, there exists a minimiser $\displaystyle \dm^*\in\argmin_{\dm\in \MC^+(\ci01\times\bar\Omega)} \C{E}(\dm)$ such that
	\begin{equation}\label{wfrexample-eq-structuremin-rho}
		\text{for some } a_i\geq 0,\ (\h^i,\g^i)\in \HG, \qquad \forall t\in\ci01,\qquad \dm^*_t = \sum_{i=1}^{\m(T+1)} a_i \h^i(t)\delta_{\g^i(t)}.
	\end{equation}
\end{thrm}
\begin{proof}
	Fix $\dm\in \MC^+(\X)$ with $\creg(\dm)<+\infty$. 
	\edit{}{The only aspect not directly covered by \cite{BrediesUnbalanced} is the existence of the minimal pair $(v,g)$. } To confirm this, note that the set
	\begin{equation}
		\enscond{\begin{array}{c}v\in \xLtwo_\dm(\X;\xR^d)\\ g\in \xLtwo_\dm(\X;\xR)\end{array}}{\int_\X \left[\alpha + \frac\beta2|v|^2 + \frac{\beta\delta^2}{2}g^2\right]\xdif\dm \leq \creg(\dm)+1,\ \partial_t\dm + \op{div}(v\dm)=g\dm}
	\end{equation}
	is bounded, hence compact in the weak topology of $\xLtwo_\dm$. There is therefore a weakly-convergent sequence converging to a point $(v,g)$ which achieves the desired infimum in \eqref{eq:WFR energy}. As the weak form of the continuity equation is preserved under weak limits in $(v,g)$, the triplet $(\dm,v,g)$ also satisfies the continuity equation.
\end{proof}

\subsection{Reformulation in the space of measures on paths}\label{sec: WFR measures on paths}
The results of Theorem~\Rref{thm: WFR sparsity} highlight the close relationship between the representations $\dm$ and $\cm$, \ie\ dynamical measures and measures on paths. We now reformulate the energy $\C{E}$ into an equivalent energy $\op{E}\colon\MC^+(\X)\to\xR\cup\{+\infty\}$ of the form in Section~\Rref{sec: generalisation}. This energy can be written
\begin{align}
	\forall\cm\in\MC^+(\HG),\qquad \op{E}(\cm) &\eqdef \frac12\sum_{j=0}^T\norm{A_j\cmt{t_j}-\data_j}_2^2 + \reg(\cm), \label{wfrexample-eq-energypathsconvex}
	\\\reg(\cm)&\eqdef \int_{\HG} \int_0^1\left[\alpha + \frac\beta2|\g'|^2 + \frac{\beta\delta^2}{2}\left(\frac{\h'}{\h}\right)^2\right]\h\xdif t\xdif\cm(\h,\g).\label{wfrexample-eq-regpathsconvex}
\end{align}
\begin{rmrk}
	Since for all $(\h,\g) \in \HG$, $\sqrt{\h}$ and $\sqrt{\h}\g$ are absolutely continuous, they are differentiable almost everywhere in $\ci01$, hence $\g$ is differentiable at a.e. $t$ such that $\h(t)>0$. Hence, the integrand in \eqref{wfrexample-eq-regpathsconvex} makes sense when regarded as 
	\begin{equation}\label{eq: WFR w}
		w(\h,\g) \eqdef \int_{\{\h>0\}}\left[\alpha + \frac\beta2\abs{\left(\frac{\sqrt{\h}\g}{\sqrt{h}}\right)'(t)}^2 \right]\h + {2\beta\delta^2}\abs{(\sqrt{\h})'(t)}^2 \xdif t.
	\end{equation}
\end{rmrk}
\edit{}{
	We will use the operator $\Theta\colon\MC(\HG)\to\MC(\X)$ introduced in Section~\Rref{sec: preliminary} to map between $\cm$ and $\dm$. The formula in \eqref{wfrexample-eq-regpathsconvex} comes from \cite[(3.9)]{BrediesUnbalanced}, but it was only shown that $\reg(\delta_\hg)=\creg(\Theta(\delta_\hg))$. We cannot hope that this holds on the whole of $\MC(\HG)$ because $\Theta$ is not a one-to-one mapping, $\MC(\HG)$ is a much larger space than $\MC(\X)$. The next lemma is needed to confirm the equivalence of minimisers between $\op{E}$ and $\C{E}$.}
\begin{lmm}\label{lmma: WFR sparsity sigma}
	Choosing $\varphi(\hg) = \frac{\alpha}{\op{E}(0)}$, the function $\op{E}$ is lower semi-continuous with compact sub-levelsets with sparse minimisers in $D=\enscond{\cm\in\MC^+(\HG)}{\int_\HG\varphi\xdif\cm\leq 1}$ (Theorem~\Rref{thm: structure of E}).
	Also, for any $\dm\in\MC^+(\X)$,
	\begin{equation}\label{eq: WFR W = min W}
		\creg(\dm) = \min_{\cm\in\MC^+(\HG)}\enscond{\reg(\cm)}{\Theta(\cm)=\dm}
	\end{equation}
	where $\min\emptyset=+\infty$. We conclude that
	\begin{equation}
		\min\enscond{\op{E}(\cm)}{\cm\in\MC^+(\HG)} = \min\enscond{\C{E}(\dm)}{\dm\in\MC^+(\X)},
	\end{equation}
	and 
	\begin{equation}
		\cm\in \argmin\enscond{\op{E}(\cm)}{\cm\in\MC^+(\HG)} \iff \Theta(\cm)\in \argmin\enscond{\C{E}(\dm)}{\dm\in\MC^+(\X)}.
	\end{equation}
\end{lmm}
\begin{proof}
	The only requirement for Theorem~\Rref{thm: structure of E} which is not explicitly assumed is that the function $w$ is lower semi-continuous with compact sub-levelsets. This is proved the appendix in Lemma~\Rref{lmma: WFR is coercive and lsc}.
	
	Secondly, fix $\dm\in \MC^+(\X)$. If $\creg(\dm)=+\infty$, then $\creg(\dm) \geq \min\enscond{\reg(\cm)}{\Theta(\cm)=\dm}$ is clear. Otherwise, let $\hat\cm\in\MC^+(\HG)$ be the measure given by Theorem~\Rref{thm: WFR sparsity} satisfying $\dm=\Theta(\hat\cm)$, then
	\begin{align}
		\creg(\dm) &= \int_\X \left[\alpha + \frac\beta2|v|^2 + \frac{\beta\delta^2}{2}g^2\right]\xdif\dm &\text{Theorem \Rref{thm: WFR sparsity}}
		\\&= \int_\HG\int_0^1 \left[\alpha + \frac\beta2|v(t,\g(t))|^2 + \frac{\beta\delta^2}{2}g(t,\g(t))^2\right]\h(t)\xdif t\xdif\hat\cm(\h,\g) &\text{Theorem \Rref{thm: Theta properties}(2)}
		\\&= \int_\HG\int_0^1 \left[\alpha + \frac\beta2|\g'(t)|^2 + \frac{\beta\delta^2}{2}\left(\frac{\h'(t)}{\h(t)}\right)^2\right]\h(t)\xdif t\xdif\hat\cm(\h,\g) = \reg(\hat\cm). &\eqref{eq: WFR v = x'}
	\end{align}
	This confirms $\creg(\dm) \geq \min\enscond{\reg(\cm)}{\Theta(\cm)=\dm}$ for all $\dm\in\MC^+(\X)$.
	
	The converse holds by Jensen's inequality. In particular, because $\creg$ is convex, proper, and lower semi-continuous\cite[Lemma A.6]{BrediesUnbalanced}, by \cite[Thm. 5]{Rockafellar1974} there exists a collection $\{(c_i,\psi_i)\}_{i\in I}\subset \xR\times\Co(\X)$ such that 
	\begin{equation}
		\forall\dm\in\MC^+(\X),\qquad \creg(\dm) = \sup_{i\in I} \left(c_i + \int_\X \psi_i\xdif\dm\right).
	\end{equation}
	As $\creg$ is positively homogeneous, also $c_i=0$. For any $\hat\cm\in\MC^+(\HG)$, $i\in I$, we have $\Theta(\hat\cm)\in\MC^+(\X)$ and
	\begin{multline}
		\int_\X \psi_i\xdif\Theta(\hat\cm) = \int_\HG\int_0^1 \h(t)\psi_i(t,\g(t))\xdif t\xdif\hat\cm
		\leq \int_\HG\sup_{j\in I}\int_0^1 \h(t)\psi_j(t,\g(t))\xdif t\xdif\hat\cm
		\\= \int_\HG\left[\sup_{j\in I}\int_\X \psi_j\xdif \Theta(\delta_\hg)\right]\xdif\hat\cm
		= \int_\HG\creg(\Theta(\delta_\hg))\xdif\hat\cm.
	\end{multline}
	It is shown in \cite[Prop. 3.9]{BrediesUnbalanced} that $\creg(\Theta(\delta_{\hg})) = w(\hg)$ from \eqref{eq: WFR w}. Now taking a supremum over $i\in I$ gives
	\begin{equation}
		\creg(\Theta(\hat\cm)) \leq \int_\HG\creg(\Theta(\delta_\hg))\xdif\hat\cm = \int_\HG w(\hg)\xdif\hat\cm = \reg(\hat\cm).
	\end{equation}
	This shows $\creg(\dm) \leq \min\enscond{\reg(\cm)}{\Theta(\cm)=\dm}$ for all $\dm\in\MC^+(\X)$, which confirms \eqref{eq: WFR W = min W} when combined with the ``$\geq$'' result. Note that $\cmt{t}=\Theta(\cm)$ for all $\cm\in\MC(\HG)$, $t\in\ci01$, so \eqref{eq: WFR W = min W} also holds for $\C{E}$ and $\op{E}$.

	Finally we consider the equivalence of minimums and minimisers. From Theorems~\Rref{thm: structure of E} and \Rref{thm: WFR sparsity} we know that there are two minimisers $\cm^*$ and $\dm^*$ respectively. From the first part of this proof there exists some $\hat\cm\in\MC^+(\HG)$ with $\dm^*=\Theta(\hat\cm)$ and $\creg(\dm^*) = \reg(\hat\cm)$. Recall that $\cmt[\cm^*]{t}=\Theta(\cm^*)_t$ and equivalently for $\hat\cm$, so \eqref{eq: WFR W = min W} also holds replacing $\reg/\creg$ with $\op{E}/\C{E}$. We conclude that
	\begin{equation}
		\inf_{\cm\in\MC^+(\HG)} \op{E}(\cm) = \op{E}(\cm^*) \geq \C{E}(\Theta(\cm^*)) \geq \inf_{\dm\in\MC^+(\X)}\C{E}(\dm) = \C{E}(\dm^*) = \op{E}(\hat\cm) \geq \inf_{\cm\in\MC^+(\HG)} \op{E}(\cm),
	\end{equation}
	so both minimal energies are equal and achieved by minimisers $\Theta(\cm^*)$ and $\hat\cm$ respectively.
\end{proof}

\subsection{Discrete-time formulation}\label{sec:WFR discrete-time energy}
Problems of the form \eqref{wfrexample-eq-energypathsconvex} also have a discrete-time structure inherited from discrete-time data. The same argument from Section~\ref{sec: discrete time formulation} is valid for the WFR penalty, although it is very hard to find the explicit form. Analytically the WFR penalty can be expressed in the form required for Assumption~\Rref{eq: step penalty}. In particular, for 
\begin{equation}\label{eq:WFR exact step}
	\step_j(\hg_{j-1},\hg_j) = \inf_{\hg \in \HG} \enscond{\int_{t_{j-1}}^{t_j} \left[\alpha + \frac\beta2|\g'(t)|^2 + \frac{\beta\delta^2}{2}\left(\frac{\h'(t)}{\h(t)}\right)^2 \right]\h(t)\xdif t 
	}{\hg = (\h,\g),
	\begin{array}{c}\hg(t_{j-1}) = \hg_{j-1},\\ \hg(t_j) = \hg_j\end{array}},
\end{equation}
geodesics of the WFR penalty satisfy
\begin{equation}
	w(\hg) = \int_0^1\left[\alpha + \frac\beta2|\g'|^2 + \frac{\beta\delta^2}2\left(\frac{\h'}{\h}\right)^2\right]\h\xdif t = \sum_{j=1}^T \step_j(\hg(t_{j-1}),\hg(t_j)).
\end{equation}
This formula can be verified with the Euler-Lagrange equation. The key point is that the left-hand side only involves up to first order derivatives so only the zeroth order constraints are needed to interpolate on intervals $\oi{t_{j-1}}{t_j}$. The exact shape or formulae for the WFR-geodesics is not so convenient as for the Benamou-Brenier penalty ($\delta\to+\infty$), although for $\alpha=0$ it is known from (\cite[Thm. 5.6]{Chizat2018})
\begin{equation}\label{eq:wfr:closedform}
	\step_j((\h_{j-1},\g_{j-1}),(\h_j,\g_j)) = \frac{4\beta\delta^2}{t_j-t_{j-1}}\left[\tfrac{\h_j+\h_{j-1}}{2} - \sqrt{\h_j\h_{j-1}}\cos\left(\min\left(\tfrac{|\g_j-\g_{j-1}|}{2\delta},\pi\right)\right)\right].
\end{equation}
More details can be found in Corollary~4.1(i) of the preprint (\url{https://arxiv.org/pdf/1506.06430v2.pdf}) of \cite{Chizat2018b}. In summary, the map $\hg\mapsto\argmin_{\tilde\hg\in\HG} \enscond{w(\tilde\hg)}{\tilde\hg(t_j)=\hg(t_j)}$ is continuous, $\tilde \h(t)$ is a simple quadratic on each $\oi{t_j}{t_{j+1}}$, and $\tilde\g(t)$ follows a straight line between $\g(t_j)$ and $\g(t_{j+1})$ with speed varying like arctan.

%% file: results.tex

\section{Numerical results}\label{sec: results}
For numerical experiments we implement variants of Algorithm~\Rref{alg: abstract FW} using the linear-oracle strategy discussed in Section~\Rref{sec: dynamical programming}. The expanded form of this algorithm is given in Algorithm~\Rref{alg: final algorithm}. The choice of $A_j$, $\op{F}_j$, $t_j$, and $\step_j$ define the optimisation problem, then the final choices of $\mesh_j$ and $k$ dictate the variant of the algorithm. We always choose the sliding step to select a local minimum of $\op{E}$ (or $\tilde{\op{E}}$) as computed by an implementation of the L-BFGS algorithm. All code to reproduce the results and figures in this work can be found online\footnote{\url{https://gitlab.inria.fr/rtovey/DP-for-dynamic-IPs}}.

The classical sliding Frank-Wolfe algorithm (Algorithm~\Rref{alg: abstract FW}) can be recovered by choosing $\mesh_j=\ci{-1}{1}\times\bar\Omega$ and $k=1$. The potential for $k>1$ was considered in \cite[Sec 5.1.5]{BrediesNumerics} as a ``multistart'' parameter. The idea is that the approximation of optimal curves (\ie\ $Y^*[\cdot,T]$) is very expensive, so some of the other near-optimal curves should also be used to improve efficiency of the algorithm. Suppose $\mesh_j$ is chosen to approximate $\ci{-1}{1}\times\bar\Omega$ with a grid of $M$ masses in $\ci{-1}{1}$ and $N^d$ points in $\bar\Omega$. Assuming $T$ is fixed, Theorem~\Rref{thm: linear oracle complexity} states that the computational complexity of the linear oracle at every iteration is $O((N^dM)^2)$. If the ``true curves'' are easy to find, then we can hope to reduce this to $O((N^dM)^2k^{-1})$ per curve (see \cite{Ding2020,flinth_linear_2021}).

We consider two families of algorithms implementing Algorithm~\Rref{alg: final algorithm}, one stochastic and the other deterministic. Throughout this section we fix $\Omega\eqdef \oi01^2$ and only consider non-negative curves $\h\geq0$, so the algorithms can be stated as: 
\begin{description}
	\item[$(k,N,M)$-random mesh] The multistart parameter is $k$. For each $n$ and $j$ we generate new independent uniformly random points $H=\{\h_i^j\sim \C U[\ci01]\}_{i=1}^M$, $X = \{\x_i^j\sim \C U[\ci01^2]\}_{i=1}^{N^2}$ and 
	\begin{equation}
		\mesh_j \eqdef \enscond{(\h,\x)}{\h\in H,\ \x\in X}.
	\end{equation}
	When $M=0$ we take $H=\{1\}$ (\ie\ a balanced mesh).

	\item[$(k,N,M)$-uniform mesh] The multistart parameter is $k$. For each $n$ we choose $\mesh_0=\ldots=\mesh_T$ to be of the form
	\begin{equation}
		\mesh_j = \mesh_j(\tilde N) \eqdef \enscond{(\h,(\x_1,\x_2))}{\h\in\{0,M^{-1},\ldots,1\},\ \x_1,\x_2\in\{0,\tilde N^{-1},\ldots,1\}}
	\end{equation}
	for some $\tilde N\geq1$. We start with $\tilde N=16$ at $n=0$, then increment $\tilde N\gets 2\tilde N$ whenever $\op{E}(\cm^{n+1}) = \op{E}(\cm^n)$ (tested after line 13 of Algorithm~\Rref{alg: final algorithm}). This process is terminated once $\tilde N>N$. Again, if $M=0$ then we only allow the balanced mass $\h=1$.
\end{description}
For all problems, $(1,+\infty,+\infty)$-random and $(1,+\infty,+\infty)$-uniform mesh algorithms are equivalent to the exact Algorithm~\Rref{alg: abstract FW}. For balanced problems, such as in Section~\Rref{sec:BB example intro}, it is sufficient to use the triplet $(1,+\infty,0)$. The random algorithms with $N,M\geq 1$ are guaranteed to converge (eventually) to the exact minimiser by Lemma~\Rref{lmma: a.s. FW convergence}. On the other hand, the uniform path algorithms have nice computational properties which allows $Y^*[\cdot,T]$ to be computed more efficiently at larger $N$. Whilst all results for random meshes are asymptotic, the uniform path algorithms also provide a clear stopping criterion. The algorithm stops at iteration $n$ when $\op{E}(\cm^{n+1})=\op{E}(\cm^n)$, often this will mean $\cm^{n+1}=\cm^n$. Let $D_d\eqdef \enscond{\cm\in D}{\op{supp}(\cm)\subset \prod_{j=0}^T\mesh_j(N)}$ be the discretisation of $D$. Expanding on the linesearch computation, due to the convexity and smoothness of $\op{E}$,
\begin{equation}
	\forall\lambda\in\ci01,\qquad \op{E}((1-\lambda)\cm^n+\lambda\ep^{n+1}) - \op{E}(\cm^n) = \lambda\int_\HG(\op{F}'(\cm^n)+w)\xdif[\ep^{n+1}-\cm^n] + O(\lambda^2).
\end{equation}
If the linesearch terminates with $\lambda=0$, then clearly $\int_\HG(\op{F}'(\cm^n)+w)\xdif[\ep^{n+1}-\cm^n]\geq0$. The optimality given by the dual gap (see \cite{Jaggi2013}), can therefore be stated
\begin{equation}
	\op{E}(\cm^n) - \min_{\cm\in D_d}\op{E}(\cm) \leq \sup_{\cm\in D_d} \int_\HG(\op{F}'(\cm^n)+w)\xdif[\cm^n-\cm] = \int_\HG(\op{F}'(\cm^n)+w)\xdif[\cm^n-\ep^{n+1}] \leq 0.
\end{equation}	
We conclude that $\cm^n$ is at least optimal up to a spatial resolution of $\frac1N$.

\begin{algorithm}
	\begin{algorithmic}[1]
		\State Set $ \cm^0 \gets 0\in\MC^+(\HG)$, $n\gets0$, fix $k\in\xN$
		\Repeat
		\State Let $\eta_j= A_j^*\nabla \op{F}_j(A_j\cmt[\cm^n]{t_j})\in\Co(\bar\Omega)$ for $j=0,\ldots,T$, and \raggedright \centerline{$\tilde{\op{E}}((\h,\g)) = \sum_{j=0}^T \eta_j(\g_j)\h_j + \sum_{j=1}^T\step_j((\h_{j-1},\g_{j-1}),(\h_j,\g_j))$}
		\State Choose $\mesh_j\subset\ci01\times\bar\Omega$, $j=0,\ldots T$ \Comment{the discrete mesh}
		\State Compute $Y^*[y,T]\in\prod_{j=0}^T\mesh_j$ for all $y\in \mesh_T$, as in \eqref{eq: inductive extreme point computation} \Comment{discrete optimal paths}
		\State Choose $\tilde\hg^1,\ldots,\tilde\hg^k\in\op{Image}(Y^*[\cdot,T])$ with least energy in $\tilde{\op{E}}$ \Comment{select best $k$ endpoints}
		\State Find $\hg^1,\ldots\hg^k\in\ci01\times\bar\Omega$ with $\tilde{\op{E}}(\hg^i) \leq \tilde{\op{E}}(\tilde \hg^i)$ \Comment{sliding step on linearised problem}
		\State Set $\cm \gets \cm^n$, re-order index $i$ such that $\tilde{\op{E}}(\hg^1)\leq \tilde{\op{E}}(\hg^2),\ldots$
		\For{$i=1,\ldots,k$}
			\State $\lambda \gets \argmin_{\lambda\in\ci01} \op{E}((1-\lambda)\cm + \lambda \varphi(\hg^i)^{-1}\delta_{\hg^i})$ \Comment{exact linesearch}
			\State $\cm \gets (1-\lambda)\cm + \lambda \varphi(\hg^i)^{-1}\delta_{\hg^i}$
		\EndFor
		\State Choose $\cm^{n+1}$ such that $\op{E}(\cm^{n+1}) \leq \op{E}(\cm)$ \Comment{sliding step on exact problem}
		\State $n\gets n+1$
		\Until{converged}
	\end{algorithmic}
	\caption{Inexact sliding Frank-Wolfe algorithm}\label{alg: final algorithm}
\end{algorithm}

\begin{figure}\centering
	\includegraphics[width=.8\linewidth]{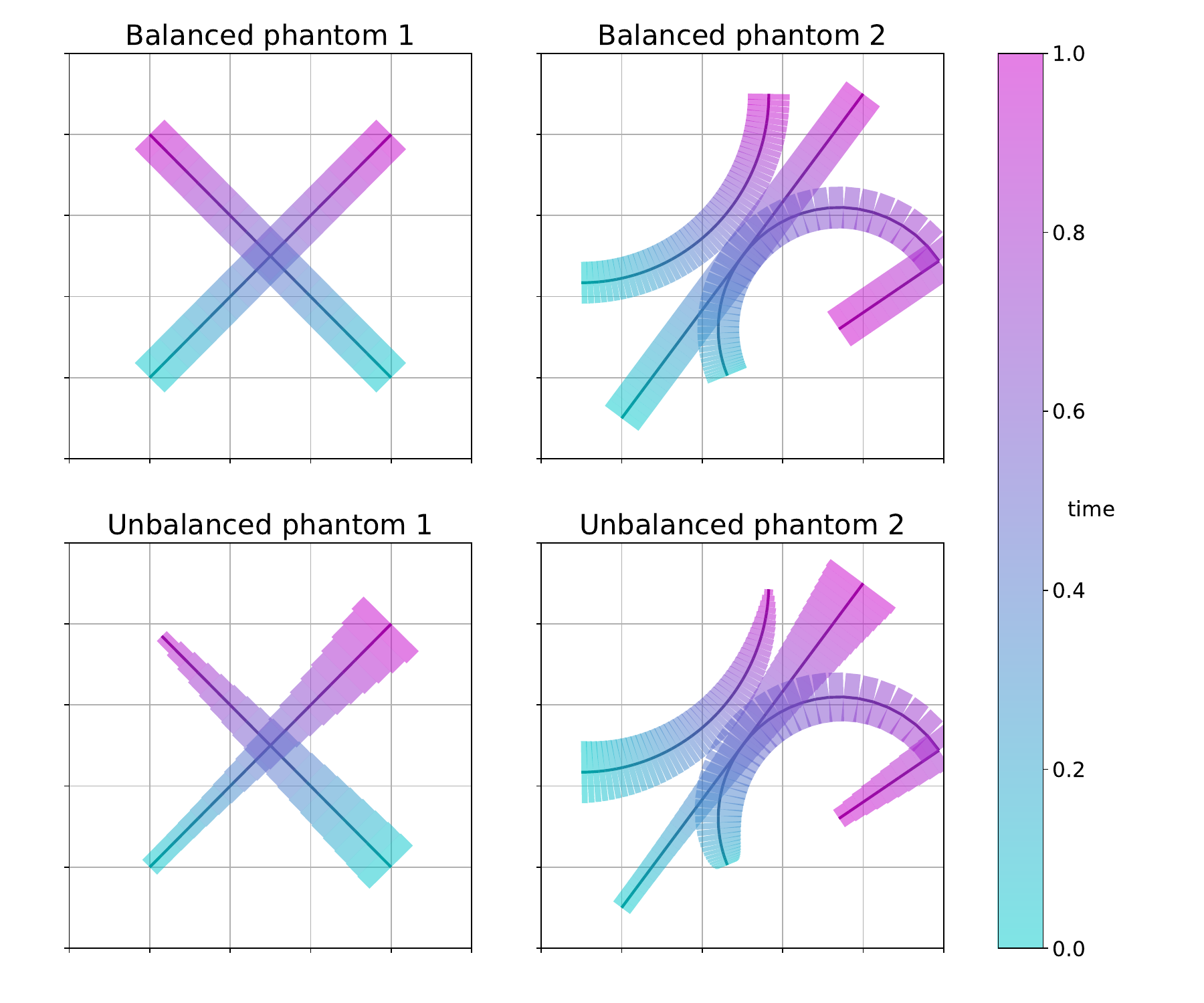}
	\caption{Synthetic phantoms of the form $\cm = \sum_i \delta_{(\h^i,\g^i)}$ for balanced (top row) and unbalanced (bottom row) examples. Colour indicates the time $t$, the solid line indicates the positions $\g^i(t)$, and the width of the overlayed band is proportional to $\h^i(t)$.}\label{fig: GT images}
\end{figure}

\subsection{Benamou-Brenier example}\label{sec: numerical BB}
First we compare directly with the numerical results presented in \cite{BrediesNumerics} for the model discussed in Section~\Rref{sec:BB example intro}. In particular, the energy we seek to minimise is $\op{E}\colon\MC^+(\HG)\to\ho{-\infty}{+\infty}$ defined by 
\begin{equation}\label{eq:numericalBB:continuous}
	\op{E}(\cm) = \frac12\sum_{j=0}^T\norm{A_j\cmt{t_j}-b_j}_2^2 + \int_\HG\int_0^1 \left(\alpha + \frac{\beta}{2}|\g'(t)|^2\right)\xdif t\xdif\cm(\h,\g)
\end{equation}
where $t_j = \frac{j}{T}$ for each $j=0,\ldots,T$,
\begin{equation}
	\HG \eqdef \enscond{(\h,\g)\in\{1\}\times\ACtwo(\ci01;\bar\Omega)}{\mbox{$\g'$ is constant on $\oi{t_{j-1}}{t_{j}}$, $1\leq j\leq T$}},
\end{equation}
and $A_j\colon \MC(\ci01^2)\to\xR^\m$ represents a finite number of smoothed Fourier samples. The precise details are given in \cite[Sec. 6]{BrediesNumerics}. In view of the convexity of $\Omega$, \eqref{eq:numericalBB:continuous} reformulates to
\begin{equation}
	\op{E}(\cm) = \frac12\sum_{j=0}^T\norm{A_j\cmt{t_j}-b_j}_2^2 + \int_\HG \left(\sum_{j=1}^{T}(t_j-t_{j-1})\left( \alpha + \frac{\beta}{2}\frac{|\g(t_j)-\g(t_{j-1})|^2}{(t_j-t_{j-1})^2}\right)\right)\xdif\cm(\h,\g).
\end{equation}

\begin{figure}\centering
	\begin{subfigure}{.7\textwidth}\centering
		\includegraphics[width=\linewidth]{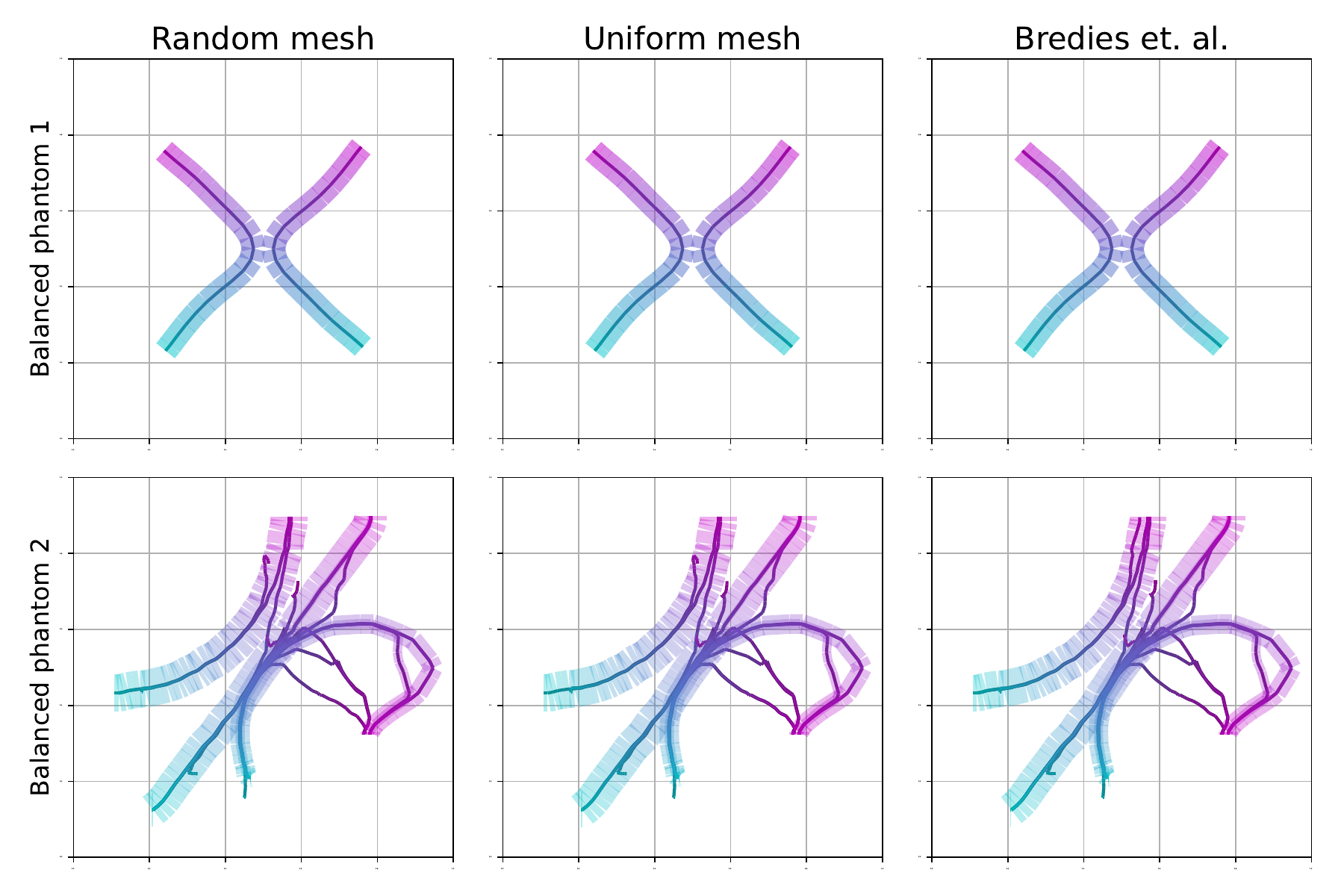}
		\caption{Final reconstructions visualised equivalently to those in Figure~\Rref{fig: GT images}.}\label{fig: balanced reconstructions}
	\end{subfigure}
	\begin{subfigure}{.7\textwidth}\centering
		\includegraphics[width=\linewidth]{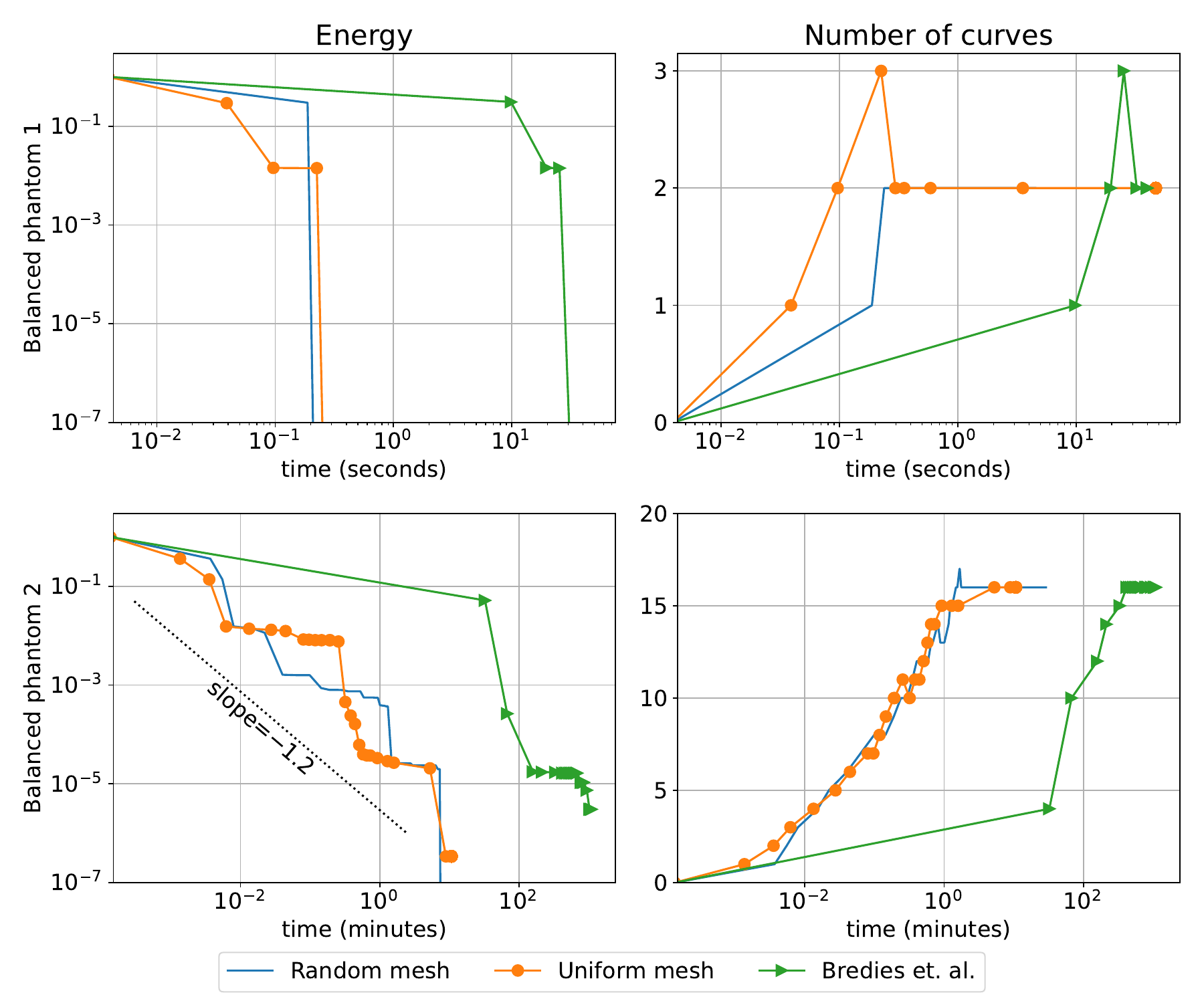}
		\caption{Convergence of energy and sparsity of reconstructions. For each phantom, every energy is translated by the smallest energy found by any method.}\label{fig: balanced convergence}
	\end{subfigure}
	\caption{Comparison of algorithm from \cite{BrediesNumerics}, the $(1,25,0)$-random mesh algorithm, and $(1,256,0)$-uniform mesh algorithm applied to balanced phantoms 1 and 2 (first and second rows respectively).}\label{fig: general balanced figure}
\end{figure}

The two phantoms are also from \cite[Sec. 6]{BrediesNumerics}. In the notation of this work, we would say that, for example, phantom 1 is represented by
\begin{align}
	\cm = \delta_{(\h^1,\g^1)} + \delta_{(\h^2,\g^2)}\qtext{where} \h^1=\h^2=1,\ \g^1(t) &= (0.2+0.6t,\ 0.2+0.6t),
	\\\g^2(t) &= (0.8-0.6t,\ 0.2+0.6t).
\end{align}
Similarly phantom 2 is the sum of 3 Dirac masses, both phantoms are shown here in Figure~\Rref{fig: GT images}. In the setting of Section~\Rref{sec: dynamical programming}, we choose $\varphi_0=0.1$ (as 10 is much larger than 2 or 3 which is the mass of phantoms 1 and 2 respectively), and 
\begin{equation}
	\step_j((\h_0,\g_0),(\h_1,\g_1)) = \alpha(t_j-t_{j-1}) + \frac\beta2\frac{|\g_1-\g_0|^2}{t_j-t_{j-1}} = \frac{\alpha}{T} + \frac{\beta T}{2}|\g_1-\g_0|^2.
\end{equation}
For phantom 1 we have $\alpha=\beta=0.5$, $T=21$, and $\alpha=\beta=0.1$, $T=51$ for phantom 2. The algorithm of \cite{BrediesNumerics} found online\footnote{\url{https://github.com/panchoop/DGCG\_algorithm/commit/553a564fd8641abcfac6067ebf51a900a6a91d0f}} is run with original parameters as a baseline, although with a time limit of 5 days when necessary. This is compared with multiple variants of the random and uniform algorithms described at the beginning of the section. The uniform algorithms are run to convergence and the random variants are run for 100 and 10,000 iterations for phantoms 1 and 2 respectively. Figure~\Rref{fig: balanced reconstructions} shows the reconstructions of each algorithm with default parameters, each image is visually equivalent. The convergence behaviour is shown in Figure~\Rref{fig: balanced convergence}. The energy plots confirm that each reconstruction has approximately the same energy, although we find that the random mesh algorithm finds the lowest energy, closely followed by the uniform mesh. Similarly, the sparsity of each final reconstruction is equal. The greatest difference between algorithms is run-time. The random and uniform mesh algorithms are over 100 times faster than that of \cite{BrediesNumerics} in both examples.

\begin{figure}\centering
	\includegraphics[width=.9\linewidth]{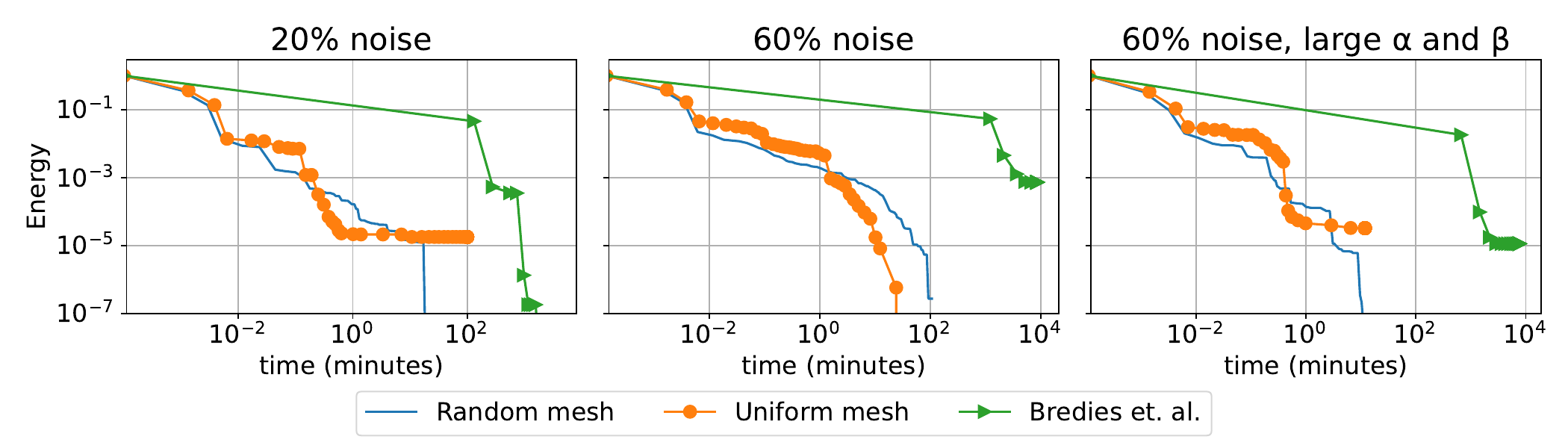}
	\caption{Convergence of three default algorithms (see Figure~\Rref{fig: general balanced figure}) with different levels of noise and choice of $\alpha=\beta\in\{0.1,0.3\}$.}\label{fig: noise convergence}
\end{figure}

We also replicate the noise-scenarios for phantom 2 as tested in \cite{BrediesNumerics}. Our only modification of the baseline algorithm is to remove the early-stopping routine and run the algorithm for the minimum of 21 iterations or 5 days (\cf \cite[Table 1]{BrediesNumerics}). We again use the $(1,25,0)$-random mesh and $(1,256,0)$-uniform mesh algorithms for comparison in each example. In the three noisy scenarios we add 20\%, 60\%, and 60\% Gaussian white noise to the data respectively. The first two scenarios use $\alpha=\beta=0.1$ while the third scenario uses $\alpha=\beta=0.3$. Seen in Figure~\Rref{fig: noise convergence}, both the random and uniform mesh algorithms converge with similar rates, while the algorithm of Bredies \etal\ is still at least 100 times slower. We see the expected behaviour that the uniform variant converges to a (possibly non-optimal) energy. As predicted by Theorem~\Rref{lmma: a.s. FW convergence}, the random variant often finds an even lower energy, despite having a much smaller value of $N$.
	
\layout{\pagebreak}\subsection{Wasserstein-Fisher-Rao example}
In this section we show numerical results for the unbalanced transport example presented in Section~\Rref{sec:WFR example}, the data fidelity is the same as in the Benamou-Brenier example. Ideally we would use the exact function $d_{\alpha,\beta,\delta}$ from \eqref{eq:WFR exact step} to define $\step_j$, but it lacks a closed-form expression, hence for computational reasons we use the approximation
\begin{equation}
	\step_j(\hg_0,\hg_1) \eqdef \alpha\frac{\h_0+\h_1}{2}(t_j-t_{j-1}) + d_{0,\beta,\delta}(\hg_0,t_{j-1}, \hg_1, t_j).
\end{equation}
where $d_{0,\beta,\delta}$ is given in~\eqref{eq:wfr:closedform}. As in Section~\Rref{sec: numerical BB} we use $\alpha=\beta=0.5$ or $0.1$ for the first and second phantom respectively, also $\varphi=\delta = 0.1$ throughout. This leads to the explicit form of $\step_j$
\begin{equation}
	\step_j(\hg_0,\hg_1) = \frac{\alpha}{T}\frac{\h_0+\h_1}{2} + \frac{\beta T}{25}\left[\frac{\h_0+\h_1}2 - \sqrt{\h_0\h_1}\cos\left(\min\left(5|\g_0-\g_1|,\pi\right)\right)\right].
\end{equation}
Numerically we re-parametrise the mass to $\tilde\h_j \eqdef \sqrt{\h_j}$ for each $j$ so that $\step_j$ is a $\xCone$ function of $\tilde\h_j$. Note that this only effects the sliding step of the Frank-Wolfe algorithm, the remaining steps are unchanged.

Our synthetic phantoms are equivalent to those in Section~\Rref{sec: numerical BB} but with modified, time-dependent masses. For example, the first phantom is now $\cm = \delta_{(\h^0,\g^0)} + \delta_{(\h^1,\g^1)}$ where
\begin{align}
	\h^0(t) &= \tfrac12(1+3t^2), &\g^0(t) &= (0.2+0.6t,\ 0.2+0.6t), 
	\\ \h^1(t) &= \tfrac32\sqrt{1-t}, &\g^1(t) &= (0.8-0.6t,\ 0.2+0.6t).
\end{align}
The transformation for the second phantom is very similar and the precise formulae can be found in the supplementary code. All curves $\h$ have been normalised so $\int_0^1\h\xdif t=1$, $\norm{\h}_\infty\leq 2$.

\begin{figure}\centering
	\begin{subfigure}{.7\textwidth}\centering
		\includegraphics[width=.8\linewidth]{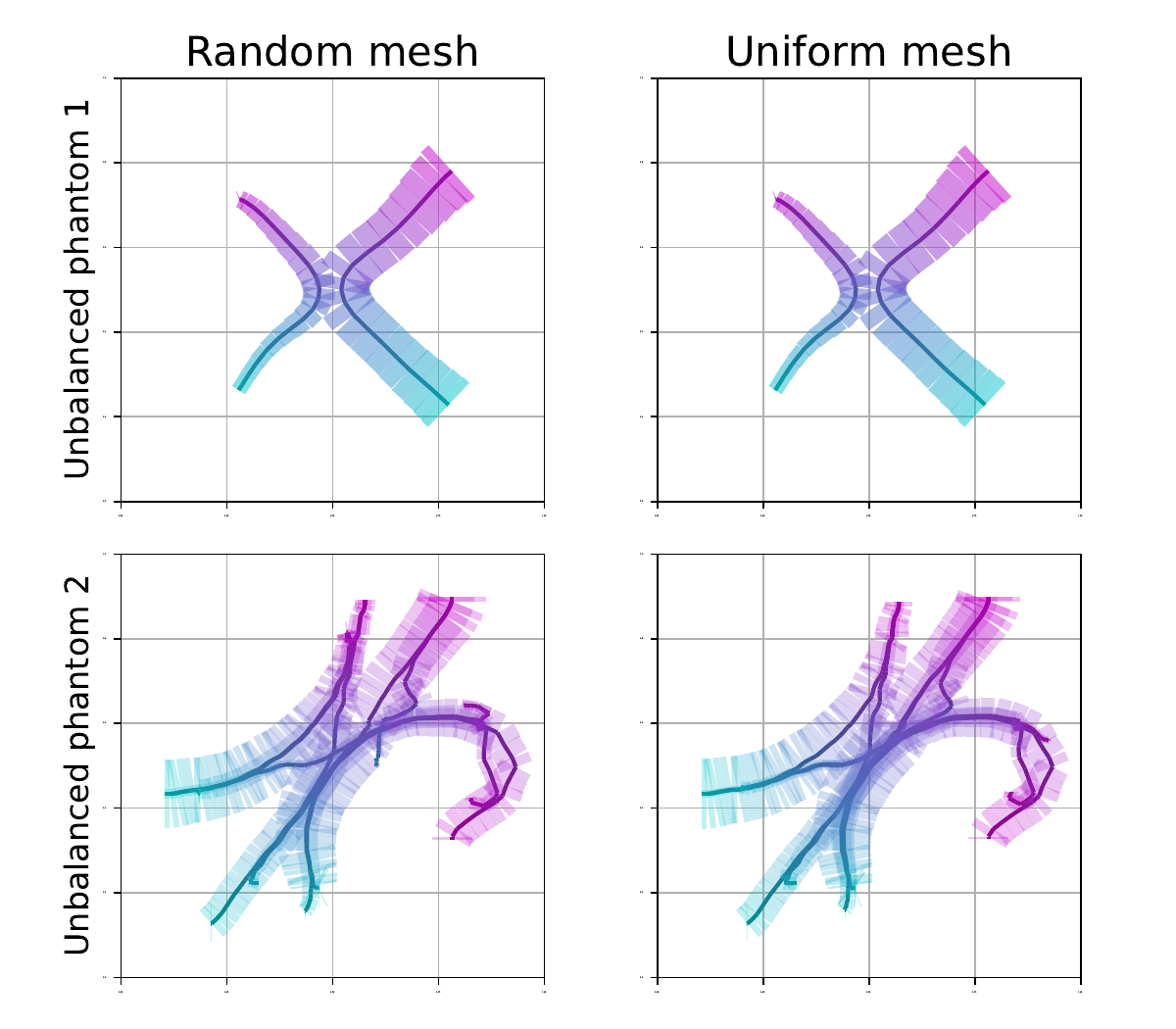}
		\caption{Final reconstructions visualised equivalently to those in Figure~\Rref{fig: GT images}.}\label{fig: unbalanced reconstructions}
	\end{subfigure}
	\begin{subfigure}{.7\textwidth}\centering
		\includegraphics[width=\linewidth]{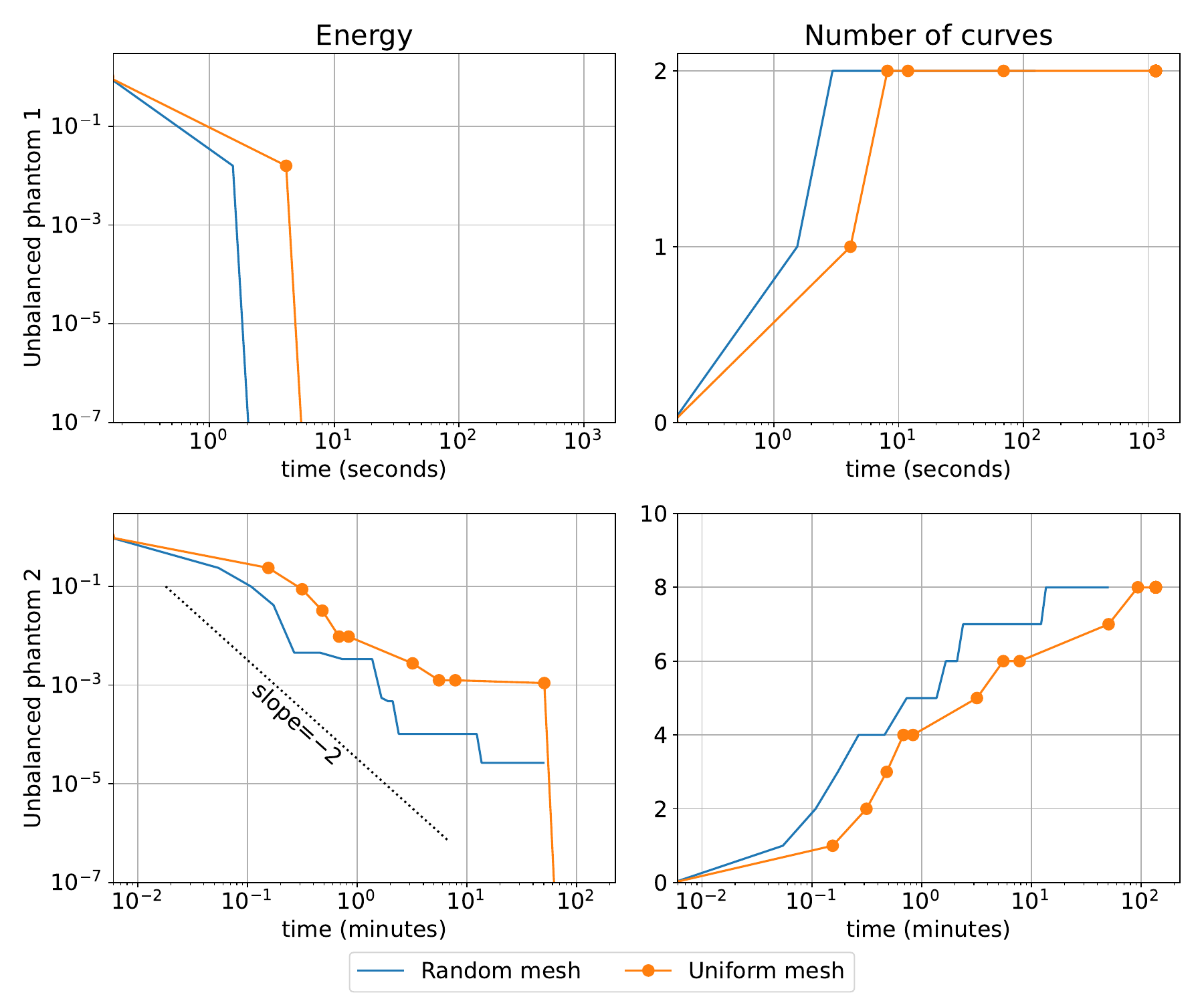}
		\caption{Convergence of energy and sparsity of reconstructions. For each phantom, the energies are translated by the smallest energy found by either method.}\label{fig: unbalanced convergence}
	\end{subfigure}
	\caption{Comparison of the $(1,25,10)$-random mesh and $(1,128,10)$-uniform mesh algorithms applied to unbalanced phantoms 1 and 2 (first and second rows respectively).}
\end{figure}

Again, we run the uniform-mesh algorithms to convergence but now the random-mesh is only run for 100 or 1,000 iterations for phantoms 1 and 2 respectively. The reconstructions are shown in Figure~\Rref{fig: unbalanced reconstructions} with corresponding convergence plots in Figure~\Rref{fig: unbalanced convergence}.

\layout{\pagebreak}\subsection{Observations on parameter choices}
Both the random and uniform mesh algorithms have three parameters to choose: the multi-start parameter $k$, spatial resolution $N$, and mass resolution $M$. The first phantom (balanced or unbalanced) nicely highlights characteristics of the reconstructions but is too trivial numerically to compare different algorithm choices, all methods converge within a few iterations. For the second phantom we prioritised the trade-off of between energy and computation time. 

The classical choice of $k=1$ was always a competitive choice. Large $k$ potentially enables the algorithm to find multiple atoms in one iteration, but slows down the sliding steps. We found that $k\in[1,5]$ was a reasonable range depending on the sparsity of the signal to be recovered.

Choice of spatial resolution made the largest impact on performance. If the resolution is too small then the algorithm will not find new atoms, but computation time of the linear oracle scales with $N^4$. For the random mesh, we found that $N=25$ was a good balance. For the uniform mesh the resolution is also a stopping criterion so we used the more conservative $N=256$ and $512$ for the balanced and unbalanced examples respectively.

In the unbalanced experiments the value of $M$ had little effect. Although we don't include this figure, we observed that the choice $M=0$ was also competitive for our unbalanced phantom 2, achieving energies within $10^{-3}\%$ of the best observed energy. It's unclear whether this will generalise to more complicated examples, we chose $M=10$ as a default. Computational complexity of the linear oracle also scales with $M^2$, so $M$ should not be too large.

Both the random and uniform mesh algorithms have advantages over each other. The main advantages of a uniform mesh are computational: one can use a finer resolution (\ie\ larger $N$), and there is a clear stopping criterion. In practice this was a very reliable method without parameter tuning and was as fast as the random algorithm. The random algorithms have analytical guarantees of converging asymptotically to the true minimiser, and indeed it achieved the best observed energy in all but one of our experiments, performing noticeably better in the noisy case. The challenge is setting a stopping criterion, in phantom 2 of Figure~\Rref{fig: unbalanced convergence} one sees several plateaux where the algorithm fails to find a descent direction for a number of iterations before continuing to descend.

\begin{figure}[H]\centering
	\includegraphics[width=.8\linewidth]{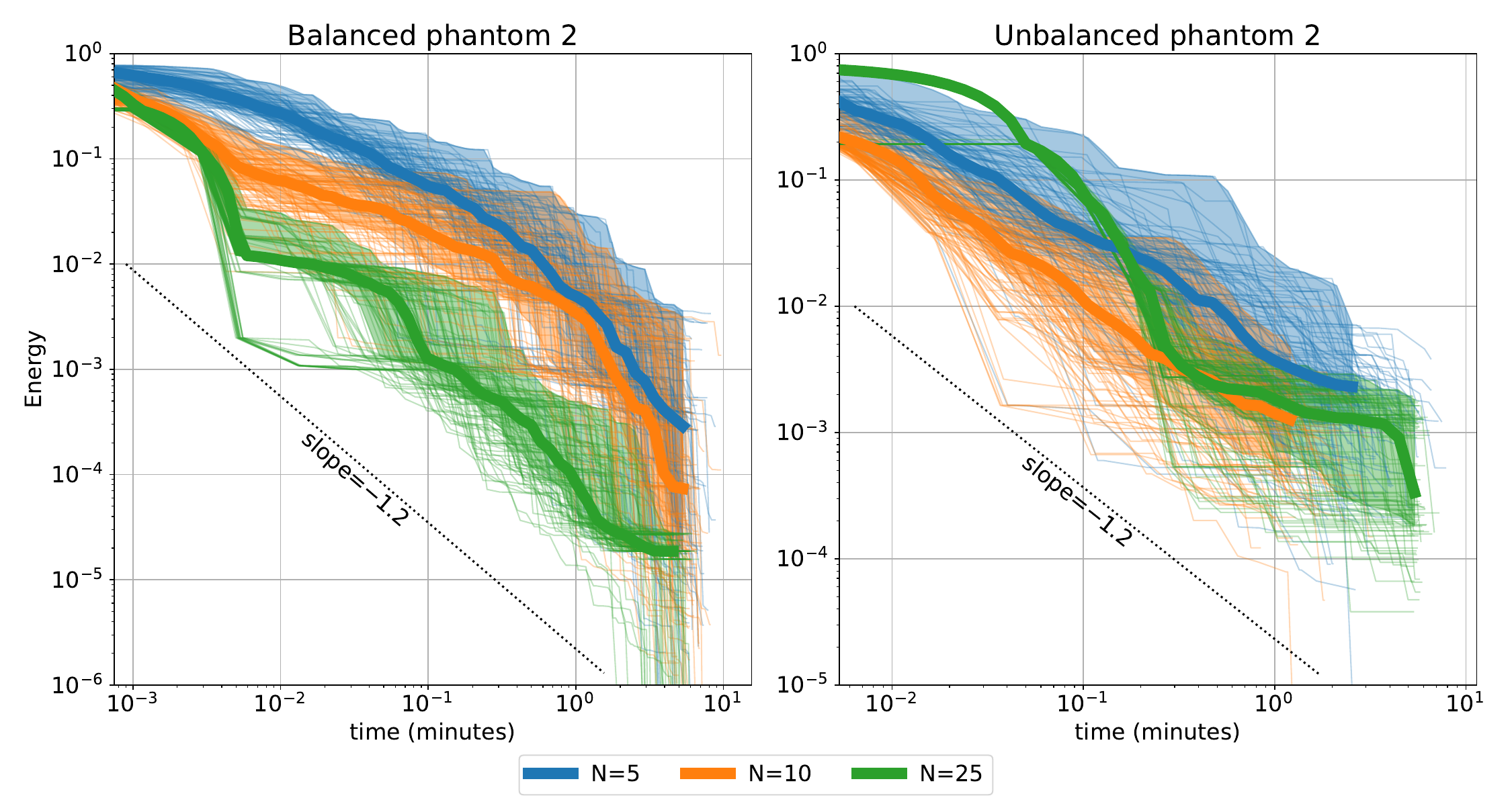}
	\caption{Convergence of different random realisations of the random mesh algorithms applied to the balanced and unbalanced phantom 2. The two algorithms are $(1,N,0)$- and $(1,N,10)$-random mesh run for 1000 and 100 iterations respectively. 100 random instances of each algorithm are run, each drawn as a thin line. The median is drawn with a thick line, and the inter-quartile range indicated by the shaded area.}\label{fig: spread convergence}
\end{figure}

Figure~\Rref{fig: spread convergence} shows the random variation of the random mesh algorithm with different values of $N\in\{5,10,15\}$. Increasing $N$ both decreases the spread and improves the performance of the algorithm. In the balanced example, the three algorithms complete 1000 iterations in the same time, showing that the linear oracle step is a small part of the total time. On the other hand, with larger $M$ the $N=25$ algorithm is nearly 10 times slower than $N=10$, indicating that the $O(N^4M^2)$ cost is starting to dominate the computation time. It's possible that $N=5$ is slower than $N=10$ for the unbalanced example because the sliding step has to work much harder to find local minima.

%% file: conclusion.tex
\section{Conclusion}\label{sec: conclusion}
The main contribution of this work was to extend a variational model proposed by Bredies \etal\ in order to accelerate the reconstruction of sparse measures in dynamical inverse problems. Using algorithms developed for computing shortest paths on graphs, we can improve the speed by a factor of 100 while still finding lower energy solutions. This allows us to process new unbalanced examples where the mass of curves is not constant in time, our proposed algorithms still recover good reconstructions in a reasonable amount of time.

We also presented new analysis of a stochastic variant of Frank-Wolfe which guarantees the convergence of our algorithm (in energy) to a globally optimal solution. This is supported by our experiments where we see that the random-mesh algorithm achieves the lowest energy in almost every example.

One feature of the algorithm in \cite{BrediesNumerics} which we did not take advantage of was the idea of importance sampling. We chose points $y\in \mesh_j$ uniformly randomly, whereas it is likely to be beneficial to choose $y$ such that $\eta_j(y)$ is small. As an example, if the step functions $\step_j$ satisfy the triangle inequality, then this bias could be implemented by choosing points $y$ such that $y$ is a local minima of
\begin{equation}
	\tilde y\mapsto \eta_j(\tilde y) + \step_{j-1}(y,\tilde y) + \step_j(\tilde y, y).
\end{equation}
In Algorithm~\Rref{alg: final algorithm}, implementing such a sliding step on $\mesh_j$ between lines 4 and 5 would preserve the analytical properties of the current algorithm, whilst possibly improving practical performance. However, it is also possible that the sliding step already present (\eg\ line 7) is powerful enough to find these improved mesh points after computing $\tilde\hg^i$, without the extra help beforehand. It is likely that the benefits depend on the application, particularly on the smoothness of the operators $A_j$.

\layout{\pagebreak}

%% file: flat_metric.tex
\section{Preliminary results}\label{app: prelims}
\subsection{Properties of the flat metric}
\begin{lmm}[Lemma~\Rref{lmma: metric properties}]\label{app: metric properties}
	Define $\metric\colon\HGF\times\HGF\to\hc{0}{+\infty}$ by
	\begin{align}
		\metric((\h_1,\g_1),(\h_2,\g_2)) &\eqdef \sup_{t\in\ci01} \op{d}_F((\h_1(t),\g_1(t)),(\h_2(t),\g_2(t))) \qquad\text{where}\\
		\op{d}_F((r_1,\x_1),(r_2,\x_2)) &\eqdef \splitln{|r_1|+|r_2|}{r_1r_2\leq 0\text{ or } |\x_1-\x_2|\geq 2}{|r_1-r_2| + \min(|r_1|,|r_2|)|\x_1-\x_2|}{\qquad\text{else,}}\label{eq:defdf}
	\end{align}
	then $(\sfrac\HGF\sim,\metric)$ is a complete separable metric space where 
	\begin{equation}
		(\h_1,\g_1)\sim (\h_2,\g_2) \iff \h_1=\h_2\qtext{and}\forall t\in\{\h_1\neq0\},\quad \g_1(t)=\g_2(t).
	\end{equation}
	Convergence of a sequence $\hg_n=(\h_n,\g_n)\in\HGF$ in the metric $\metric$ can equivalently be stated as:
	\begin{equation}\label{eq: convergence in metric}
		\left[\hg_n\stackrel{\metric}{\to}(\h,\g)\right] \iff \Big[\h_n\to\h \text{ in } \Co(\ci01)\text{ and for all $\epsilon>0$, } \g_n\to\g \text{ in } \Co(\{|\h|\geq\epsilon\})\Big]
	\end{equation}
	Furthermore, for any $\psi\in\Co(\X)$, we have $\Psi\in\Co(\ci01\times\HGF)$ where
	\begin{equation}\label{eq:def: Psi}
		\forall t\in\ci01,\ (\h,\g)\in \HGF, \qquad \Psi(t,\h,\g)\eqdef \h(t)\psi(t,\g(t)).
	\end{equation}
\end{lmm}
\begin{proof}
	\emph{Complete metric space.}
	\edit[2]{%
		To see that $(\sfrac\HGF\sim,\metric)$ is a metric space, it is sufficient to show that $(\sfrac{\xR\times\bar\Omega}{\sim},\op{d}_F)$ is a metric space, where for all $(r_1,\x_1),(r_2,\x_2) \in \xR\times\bar\Omega$,
	\begin{equation}
		(r_1,\x_1)\sim (r_2,\x_2) \iff r_1=r_2\qtext{and} \text{if $r_1\neq 0$,}\quad \x_1=\x_2.
	\end{equation}
	Following the argument of \cite[Lemma 3.3]{BrediesUnbalanced} which only considered the case $\h\geq0$, observe that for all $r_1,r_2\in\xR$, $\x_1,\x_2\in\bar\Omega$
	\begin{equation}
		\op{d}_F((r_1,\x_1),(r_2,\x_2)) = \sup_{c_1,c_2\in\xR} \enscond{r_1c_1-r_2c_2}{|c_1|,|c_2|\leq 1,\ |c_1-c_2|\leq |\x_1-\x_2|}.
	\end{equation}
	For $\min(r_1,r_2)\geq0$, this has already been shown in \cite{BrediesUnbalanced}. The case $\max(r_1,r_2)\leq 0$ is equivalent, therefore it is true for all $r_1r_2\geq 0$. For the final case, without loss of generality $r_1\geq 0\geq r_2$, so the supremum is achieved by $c_1=c_2=1$. Note from the definition that for all $r_i\in\xR$, $\x_i\in\bar\Omega$,
	\begin{equation}
		|r_1-r_2| \leq \op{d}_F((r_1,\x_1),(r_2,\x_2)) \leq |r_1|+|r_2|.
	\end{equation}
	
	It is clear that $\op{d}_F\geq0$. We must also check $\op{d}_F((r_1,\x_1),(r_2,\x_2))=0\iff (r_1,\x_1)\sim(r_2,\x_2)$. This also follows from \eqref{eq: d_F equivalence} and a case analysis of the definition of $\op{d}_F$.		
	
	Finally we verify the triangle inequality. Let $(r_1,\x_1),(r_2,\x_2),(r_3,\x_3)\in\xR\times\bar\Omega$. Let $c_1,c_3 \in \ci{-1}{1}$ such that $\abs{c_1-c_3}\leq \abs{x_1-x_3}$. By symmetry of $\op{d}_F$, we may assume without loss of generality that $c_1\leq c_3$. Then
	\begin{align}
		c_3-c_1 \leq \abs{x_1-x_3}\leq \abs{x_2-x_1}+\abs{x_3-x_2},
	\end{align}
	hence there exists $c_2 \in \ci{c_1}{c_3}$ such that $\abs{c_2-c_1}\leq \abs{x_2-x_1}$ and $\abs{c_2-c_3}\leq \abs{x_2-x_3}$. Thus,
	\begin{align}
		r_1c_1-r_3c_3= (r_1c_1-r_2c_2) + (r_2c_2 -r_3c_3) \leq \op{d}_F((r_1,\x_1),(r_2,\x_2))+\op{d}_F((r_2,\x_2),(r_3,\x_3)).
	\end{align}
	Taking the supremum over $c_1, c_3$ yields
	\begin{align}
		\op{d}_F((r_1,\x_1),(r_3,\x_3))\leq \op{d}_F((r_1,\x_1),(r_2,\x_2))+\op{d}_F((r_2,\x_2),(r_3,\x_3).
	\end{align}
	Hence $\op{d}_F$ satisfies the triangle inequality, so we conclude that both $\op{d}_F$ and $\metric$ are metrics.
	}{%
	In \cite[Prop. 3.6]{BrediesUnbalanced} it was shown that $\enscond{t\mapsto\h(t)\delta_{\g(t)}}{(\h,\g)\in\HGF,\ \h\geq0}$ is a complete separable metric space with respect to the flat metric, which is simply
	$\op{d}(\h_1\delta_{\g_1},\h_2\delta_{\g_2}) = \op{d}_\HG((\h_1,\g_1),(\h_2, \g_2))$. Reducing to $\sfrac\HGF\sim$ reduces to a single representative all couples $(\h,\g)$ which map to the same element $\h\delta_{\g}$. After removing this redundancy from $\HGF$, it is clear that $(\enscond{(\h,\g)\in\sfrac\HGF\sim}{\h\geq0}, \op{d}_\HG)$ is a complete separable metric space isometrically equivalent to that in \cite[Prop. 3.6]{BrediesUnbalanced}. In the signed case, consider $\h_i^\pm = \max(0, \pm\h_i)$, our choice of $\op{d}_\HG$ is such that for all $(\h_1,\g_1),(\h_2,\g_2)\in \HGF$
	\begin{equation}
		\op{d}_\HG((\h_1,\g_1),(\h_2, \g_2)) = \op{d}_\HG((\h_1^+,\g_1),(\h_2^+, \g_2)) + \op{d}_\HG((\h_1^-,\g_1),(\h_2^-, \g_2)).
	\end{equation}
	It follows that $(\sfrac\HGF\sim,\op{d}_\HG)$ is also a complete separable metric space.
	}
	
	\emph{Convergence.}
	First, note from the definition that for all $r_i\in\xR$, $\x_i\in\bar\Omega$,
	\begin{equation}\label{eq: d_F equivalence}
		|r_1-r_2| \leq \op{d}_F((r_1,\x_1),(r_2,\x_2)) \leq |r_1|+|r_2|.
	\end{equation}
	Now, fix a sequence $\hg_n=(\h_n,\g_n)\in\HGF$ and point $\hg=(\h,\g)\in\HGF$. Suppose $\metric(\hg_n,\hg)\to 0$, then from \eqref{eq: d_F equivalence} we have $\norm{\h_n-\h}_\infty\to0$. Also, for any $\epsilon>0$ choose $N_\epsilon\in\xN$ such that 
	\begin{equation}
		\forall n\geq N_\epsilon,\qquad \norm{\h_n-\h}_\infty\leq \frac\epsilon2 \qtext{and} \metric(\hg_n,\hg)\leq \frac\epsilon2.
	\end{equation}
	Observe that for all $t\in\{|\h|\geq\epsilon\}$ and $n\geq N_\epsilon$, if $|\g_n(t)-\g(t)|\geq 2$, then $\metric(\hg_n,\hg)\geq |\h_n(t)|+|\h(t)| > \epsilon$, contradicting the choice of $N_\epsilon$. Therefore we have the uniform bound
	\begin{equation}\label{eq: g L^infty convergence}
		\sup_{|\h(t)|\geq\epsilon}|\g_n(t)-\g(t)| \leq \sup_{|\h(t)|\geq\epsilon}\frac{2}{\epsilon}\min(|\h_n(t)|,|\h(t)|)|\g_n(t)-\g(t)| \leq \frac2\epsilon\metric(\hg_n,\hg)\stackrel{n\to+\infty}{\longrightarrow}0.
	\end{equation}
	This concludes the ``$\implies$'' direction of \eqref{eq: convergence in metric}. Conversely, suppose $\norm{\h_n-\h}_\infty\to0$ and for any $\epsilon>0$, $\norm{\g_n-\g}_{\xLinfty(\{|\h|\geq\epsilon\})}\to 0$. Fix $\epsilon\in\oi02$ and $N_\epsilon\in \xN$ such that for all 
	\begin{equation}
		\forall n\geq N_\epsilon,\qquad \norm{\h_n-\h}_\infty\leq \frac\epsilon2\qtext{and} \norm{\g_n-\g}_{\xLinfty(\{|\h|\geq\epsilon\})}\leq \epsilon.
	\end{equation}
	Then, from \eqref{eq: d_F equivalence}, for all $n\geq N_\epsilon$ we have
	\begin{align}
		\metric(\hg_n,\hg) &\leq \sup_{t\in\ci01} \left\{\begin{array}{l} |\h_n(t)|+|\h(t)| \\ |\h_{n}(t)-\h(t)| + \min(|\h_n(t)|,|\h(t)|)|\g_{n}(t)-\g(t)| \end{array}\right. & \begin{array}{c}\text{if }|\h(t)|<\epsilon, \\ \text{else} \end{array}
		\\&\leq \sup_{t\in\ci01} \left\{\begin{array}{l} 5\epsilon/2 \\ \epsilon/2+\norm{\h}_\infty\epsilon \end{array}\right. & \begin{array}{c}\text{if }|\h(t)|<\epsilon, \\ \text{else.} \end{array} &\hspace{20pt}
	\end{align}
	In either case we have $\limsup_{n\to+\infty} \metric(\hg_n,\hg)\leq O(\epsilon)$, therefore $\metric(\hg_n,\hg)\to0$ as required.	
	
	\emph{Continuity.} Fix $\psi\in\Co(\X)$ and define $\Psi$ as in \eqref{eq:def: Psi}. As we have now confirmed that $(\sfrac{\HGF}{\sim},\metric)$ is a metric space, we can use the sequential definitions of continuity. Suppose $\tau_n\in\ci01$, $\hg_n=(\h_n,\g_n)\in\HGF$ and $\tau_n\to t$, $\hg_n\stackrel{\metric}{\to}\hg$. Observe for each $n$ we have
	\begin{align}
		|\Psi(\tau_n,\hg_n)-\Psi(t,\hg)| &= |\h_n(\tau_n)\psi(\tau_n,\g_n(\tau_n)) - \h(t)\psi(t,\g(t))|
		\\&= |(\h_n(\tau_n)-\h(t)+\h(t))\psi(\tau_n,\g_n(\tau_n)) - \h(t)\psi(t,\g(t))|
		\\&\leq |\h_n(\tau_n)-\h(t)|\norm{\psi}_\infty + |\h(t)||\psi(\tau_n,\g_n(\tau_n))-\psi(t,\g(t))|.
		\\&\leq [\norm{\h_n-\h}_\infty+|\h(\tau_n)-\h(t)|]\norm{\psi}_\infty + |\h(t)||\psi(\tau_n,\g_n(\tau_n))-\psi(t,\g(t))|.
	\end{align}
	Now consider the limit of $n\to+\infty$. As $\h\in\Co(\ci01)$, $\h(\tau_n)\to\h(t)$. The characterisation of convergence in $\metric$ in \eqref{eq: convergence in metric} also confirms $\norm{\h_n-\h}_\infty\to0$ and either $\h(t)=0$, or
	\begin{equation}
		|\g_n(\tau_n)-\g(t)| \leq |\g_n(\tau_n)-\g(\tau_n)| + |\g(\tau_n)-\g(t)| \to 0.
	\end{equation}
	In either case, we see that
	\begin{equation}
		\limsup_{n\to+\infty}|\Psi(\tau_n,\hg_n)-\Psi(t,\hg)| \leq \limsup_{n\to+\infty}|\h(t)||\psi(\tau_n,\g_n(\tau_n))-\psi(t,\g(t))| = 0,
	\end{equation}
	therefore $\Psi$ is continuous at $(t,\hg)$.

\end{proof}

%% file: Theta_thm.tex
\subsection{Projection properties}
\begin{thrm}[Theorem~\Rref{thm: Theta properties}]\label{app: Theta properties}
	Let $\cm\in\MC(\HGF)$. If $\int_\HGF\norm\h_1\xdif|\cm|(\h,\g)<+\infty$, then there is a unique finite Borel measure $\Theta(\cm)\in\MC(\X)$ such that
	\begin{equation}\label{eq:app: Theta def.}
		\forall\psi\in\Co(\X),\qquad \int_\X\psi(t,\x)\xdif\Theta(\cm)(t,\x) = \int_\HGF \left(\int_0^1\h(t)\psi(t,\g(t))\xdif t\right)\xdif\cm(\h,\g).
	\end{equation}
	Moreover, 
	\begin{enumerate}
		\item The mapping $\Theta\colon \enscond{\cm\in\MC(\HGF)}{\int_\HGF\norm\h_1\xdif|\cm|<+\infty} \to\MC(\X)$ is linear.
		\item Equality~\eqref{eq:app: Theta def.} holds for all $\psi\in \xLone_{|\Theta(\cm)|}(\X)$.
		\item If $\int_\HGF\norm\h_\infty\xdif\abs{\cm}<+\infty$, then $\Theta(\cm)\in\CwX$.
		\item Suppose $\h,\g\in \ACtwo(\ci01)$ for $\cm$-a.e. $(\h,\g)\in\HGF$. If there exist Borel measurable functions $v\colon \X\rightarrow \xR^d$ and $g:\X\rightarrow \xR$ such that 
		\begin{align}
			&\h'(t) = g(t,\g(t))\h(t) \text{ for $\cm$-a.e. $(\h,\g)$ and a.e. $t\in\oi01$}, \label{eq:app: continuity eq on curve1}
			\\&\g'(t) = v(t,\g(t)) \text{ for $\cm$-a.e. $(\h,\g)$ and a.e. $t$ such that $\h(t)\neq0$},
			\\ \text{and }&\int_\HGF\int_0^1 (1+|v(t,\g(t))| + |g(t,\g(t))|)\,|\h(t)|\,\xdif t\xdif\abs{\cm}(\h,\g) <+\infty,\label{eq:app: continuity eq on curve3}
		\end{align}
		then $\int_\HGF\norm\h_\infty\xdif\abs\cm<+\infty$ and $\Theta(\cm)$ satisfies the continuity equation \eqref{eq: continuity equation}. 
		
		Conversely, given $\dm\in\MC(\X)$, if $\dm\geq0$ satisfies the continuity equation \eqref{eq: continuity equation} and 
		\begin{equation}
			\int_\X (1+|v(t,\x)|^2+|g(t,\x)|^2)\xdif\dm(t,\x) <+\infty,
		\end{equation}
		then $\dm=\Theta(\cm)$ for some $\cm\in\MC^+(\HGF)$ such that \eqref{eq:app: continuity eq on curve1}-\eqref{eq:app: continuity eq on curve3} hold and $\int_\HGF\norm\h_\infty\xdif\cm<+\infty$.
	\end{enumerate}
\end{thrm}
\begin{proof}
	By Lemma~\Rref{app: metric properties}, for any $\psi\in\Co(\X)$, the map 
	\begin{equation}
		(\h,\g)\mapsto \int_0^1\h(t)\psi(t,\g(t))\xdif t = \int_0^1\Psi(t,\h,\g)\xdif t
	\end{equation}
is continuous (hence Borel) in $\HGF$ and dominated by $\norm{\h}_1\norm{\psi}_\infty$, therefore the right-hand side of \eqref{eq:app: Theta def.} is well-defined. This induces a linear form on $\Co(\X)$ which is moreover bounded, since
	\begin{equation}
		\left|\int_\HGF\int_0^1\h(t)\psi(t,\g(t))\xdif t\xdif\cm(\h,\g)\right| \leq \norm{\psi}_\infty\int_\HGF\norm{\h}_1\xdif\abs\cm.
	\end{equation}
	By the Riesz representation theorem, that linear form is represented by a unique Radon measure $\Theta(\cm) \in \MC(\X)$. This confirms that $\Theta$ satisfies the required properties for point (1).
	
	\edit{}{Points (2)-(4) have been proved in \cite{BrediesUnbalanced} under the additional assumptions $\cm\geq0$ and $\inf_{t\in\ci01}\h(t)\geq 0$. To apply these results, we need the following modified Hahn-Jordan decomposition.}
	\begin{quote}
		\layout{\vspace{-3pt}}
		\begin{clm}
			For any $\cm\in\MC(\HGF)$ with $\int_\HGF\norm\h_1\xdif\abs\cm<+\infty$ there exists $\cm^+,\cm^-\in\MC^+(\HGF)$ such that
			\begin{equation}\label{eq-cmpm-zeromass}
				\cm^\pm\left(\enscond{(\h,\g)\in\HGF}{\inf_{t \in \ci{0}{1}}\h(t)<0}\right)=0
			\end{equation}
			with $\int_{\HGF}\norm\h_1\xdif\cm^\pm<+\infty$ and $\Theta(\cm) = \Theta(\cm^+)-\Theta(\cm^-)$.
		\end{clm}
	\begin{proof}[Proof of claim.]
		The Hahn-Jordan decomposition gives $\cm = \max(0,\cm) - \max(0,-\cm)$, therefore it is sufficient to consider the case $\cm\geq 0$. Define the maps $T^\pm\colon \HGF \rightarrow \HGF$ by
		\begin{equation}
			\forall (\h,\g) \in \HGF,\qquad	T^\pm(\h,\g) \eqdef (\max(0,\pm\h), \g).
		\end{equation}
		Because $\metric(T^\pm(\hg_1),T^\pm(\hg_2))\leq \metric(\hg_1,\hg_2)$, $T^\pm$ are continuous (therefore Borel), and we can define the image measures $\sigma^\pm\eqdef (T^\pm)_\sharp \sigma$. Since the push-forward operation does not increase the total variation, we have $\cm^+, \cm^- \in \MC^+(\HGF)$. Moreover, \eqref{eq-cmpm-zeromass} holds, and by construction of the image measure
		\begin{equation}
			\int_{\HGF} \phi(\h,\g)\xdif\cm^\pm = \int_\HGF \phi(\max(0,\pm\h),\g)\xdif\cm
		\end{equation}
		for all $\phi\in\Cb(\HGF)$, and by monotone convergence
		\begin{equation}
			\int_{\HGF}\norm\h_1\xdif\cm^\pm = \int_\HGF \norm{\max(0,\pm\h)}_1\xdif\cm \leq \int_\HGF \norm{\h}_1\xdif\cm.
		\end{equation}
		Finally, we confirm that for all $\psi \in \Co(\X)$,
		\begin{align}
			\int_\X \psi(t,\x)\xdif \Theta(\cm) &= \int_\HGF \int_0^1\h(t)\psi(t,\g(t))\xdif t\xdif\cm
			\\&= \int_\HGF \int_0^1[\max(0,\h(t))-\max(0,-\h(t))]\psi(t,\g(t))\xdif t\xdif\cm
			\\&= \int_{\HGF}\int_0^1\h(t)\psi(t,\g(t))\xdif t\xdif (\cm^+-\cm^-)
			\\&= \int_\X \psi(t,\x)\xdif (\Theta(\cm^+)-\Theta(\cm^-))
		\end{align}
		as required.
	\end{proof}\end{quote}

	\noindent \edit{}{With this decomposition, we can apply the results of \cite{BrediesUnbalanced} to each $\cm^\pm$. Point (2) becomes a direct consequence of \cite[Lemma 4.4]{BrediesUnbalanced}. Under the assumptions of point (3), the same lemma also guarantees the existence of a disintegration $\dm_t^\pm$ of $\Theta(\cm^\pm)$, in the sense of \eqref{eq: disintegration equation}.} In particular,
	\begin{equation}
		\forall\psi\in\Co(\bar\Omega),\ t\in\ci01,\qquad \int_{\bar\Omega}\psi(\x)\xdif\dm_t^\pm(\x) = \int_{\HGF}\h(t)\psi(\g(t))\xdif\cm^\pm(\h,\g).
	\end{equation}
	The same property holds for $\Theta(\cm)$ by linearity. Point (3) requires $t\mapsto\int_{\bar\Omega}\psi\xdif(\dm_t^+-\dm^-_t)$ to be continuous for all $\psi\in\Co(\bar\Omega)$. By Lemma~\Rref{app: metric properties}, for each $t\in\ci01$ the function $\Psi(t,\h,\g) \eqdef \h(t)\psi(\g(t))$ is continuous on $\ci01\times\HGF$ and 
	\begin{equation}
		\left|\int_{\bar\Omega} \psi\xdif(\dm_\tau-\dm_t)\right| \leq \int_\HGF |\h(\tau)\psi(\g(\tau)) - \h(t)\psi(\g(t))| \xdif\abs{\cm}(\h,\g) = \int_\HGF |\Psi(\tau,\hg) - \Psi(t,\hg)| \xdif\abs{\cm}(\hg).
	\end{equation}
	The integrand is pointwise bounded by $2\norm{\h}_\infty\norm{\psi}_\infty$, therefore we conclude that the limit of the integral as $\tau\to t$ is 0 by dominated convergence.
	
	\edit{}{Point (4) and its converse are addressed by \cite[Thm. 4.2]{BrediesUnbalanced}}. For the forward direction, we again consider the modified Hahn-Jordan decomposition. The assumptions \eqref{eq:app: continuity eq on curve1}-\eqref{eq:app: continuity eq on curve3} are only assumed to hold for almost every $t$ and $\h$, therefore they also hold for the choice of measures $\cm^\pm$ given in the claim replacing $\h$ with $\max(0,\pm\h)$. \edit{}{We can then apply \cite[Thm. 4.2]{BrediesUnbalanced} to each component and sum for the result. The converse is exactly the statement of \cite[Thm. 4.2]{BrediesUnbalanced}} because we assume $\dm\geq0$.
\end{proof}

\begin{lmm}[Lemma~\Rref{lmma: Gamma_infty vs Gamma_1}]\label{app: Gamma_infty vs Gamma_1}
	For each $p\in\ci{1}{+\infty}$ define the set
	\begin{equation}
		\HG_p \eqdef \enscond{\hg=(\h,\g)\in\HGF}{\norm{\h}_p\leq 1},
	\end{equation}
	then
	\begin{equation}\label{eq: Gamma_p equivalence}
		\enscond{\Theta(\cm)}{\cm\in\MC(\HGF),\ \int_{\HGF}\norm\h_p\xdif\abs\cm <+\infty} = \enscond{\Theta(\hat\cm)}{\hat\cm\in\MC(\HG_p)}
	\end{equation}
	and $\Theta\colon\MC(\HG_p)\to\MC(\X)$ is narrowly continuous. 
	
	Furthermore, if $p=+\infty$, then $\forall t\in\ci01$, $\cmt[\colon]{t}\MC(\HG_\infty)\to\MC(\bar\Omega)$ is also narrowly continuous. 
	
	In particular, sequentially we have that, for any sequence $\cm^n\xrightharpoonup{*}\cm$ narrowly in $\MC(\HG_p)$:
	\begin{align}
		\text{for all }p\in\ci{1}{+\infty}, &\qquad \Theta(\cm^n)\stackrel{*}{\rightharpoonup}\Theta(\cm) \text{ narrowly in }\MC(\X),
		\\\text{if }p=+\infty,\ \forall t\in\ci01,&\qquad \cmt[\cm^n]{t}\stackrel{*}{\rightharpoonup}\cmt{t} \text{ narrowly in }\MC(\bar\Omega).
	\end{align}
\end{lmm}
\begin{proof}
	The ``$\supset$'' inclusion is clear: if $\hat\cm\in\MC(\HG_p)$, extend $\hat\cm$ by 0 such that $\hat\cm\in\MC(\HGF)$, then 
	\begin{equation}
		\int_{\HGF}\norm\h_p\xdif\abs{\hat\cm} \leq \int_{\HG_p}1\xdif\abs{\hat\cm} = \norm{\hat\cm}<+\infty
	\end{equation}
	as required. Conversely, suppose $\int_{\HGF}\norm\h_p\xdif\abs{\cm} <+\infty$.
	Note as $\h\mapsto\max(1,\norm\h_p)$ and $\h\mapsto\frac{\h}{\max(1,\norm\h_p)}$ are continuous in $\Co(\ci01)$ and $\metric(\hg_1,\hg_2)\geq \norm{\h_1-\h_2}_\infty$, both $(\h,\g)\mapsto \max(1,\norm\h_p)$ and $(\h,\g)\mapsto (\frac{\h}{\max(1,\norm\h_p)},\g)$ are continuous w.r.t. $\metric$. Define $\hat\cm$ through the (rescaled) push-forward $\hat\cm(\h,\g) = \max(1,\norm\h_p)\cdot T_\sharp\cm(\h,\g)$ where $T\colon\HGF\to\HG_p$ is defined by the map $T(\h,\g) \eqdef \left(\frac{\h}{\max(1,\norm\h_p)},\g\right)$. $T$ is Borel measurable, therefore $\hat\cm$ is a Borel measure which satisfies
	\begin{equation}
		\forall\phi\in\Cb(\HG_p),\qquad \int_{\HG_p} \phi\xdif\hat\cm = \int_{\HGF}\max(1,\norm{\h}_p)\phi\left(\frac{\h}{\max(1,\norm{\h}_p)},\g\right)\xdif \cm(\h,\g).
	\end{equation}
	Furthermore, $\hat\cm\in\MC(\HG_p)$ as 
	\begin{equation}
		\left|\int_{\HG_p} \phi\xdif\hat\cm\right| \leq \norm{\phi}_\infty\int_{\HGF}(1+\norm{\h}_p)\xdif\abs\cm(\h,\g),
	\end{equation}
	so $\norm{\hat\cm}<+\infty$.	The last step is to confirm $\Theta(\hat\cm)=\Theta(\cm)$. For any $\psi\in\Co(\X)$ observe
	\begin{align}
		\int_\X \psi(t,\x)\xdif\Theta(\hat\cm)(t,\x) &= \int_{\HG_p}\left(\int_0^1\h(t)\psi(t,\g(t))\xdif t\right)\xdif\hat\cm(\h,\g)
		\\&= \int_{\HGF} \max(1,\norm{\h}_p) \left(\int_0^1\frac{\h(t)}{\max(1,\norm{\h}_p)} \psi(t,\g(t))\xdif t\right)\xdif\cm(\h,\g) 
		\\&= \int_\X \psi(t,\x)\xdif\Theta(\cm)(t,\x).
	\end{align}
	This shows the ``$\subset$'' direction, therefore we conclude the equality in \eqref{eq: Gamma_p equivalence}. For continuity of $\Theta$, \cite[Thm. 1, Sec. 1.6]{Yosida1980} states that $\Theta$ is continuous in $(\Co(\X))'$ if and only if for every $\psi\in\Co(\X)$ there exists a finite collection $\Phi\subset\Cb(\HG_p)$ such that 
	\begin{equation}
		\forall\cm\in\MC(\HG_p),\qquad \left| \int_\X \psi(t,\x)\xdif\Theta(\cm)\right| \leq \max_{\phi\in\Phi}\left|\int_{\HG_p}\phi(\h,\g)\xdif\cm(\h,\g)\right|.
	\end{equation}
	From Lemma~\Rref{lmma: metric properties}, define $\Psi(t,\h,\g) = \h(t)\psi(t,\g(t))$, then $\phi(\h,\g) \eqdef \int_0^1\Psi(t,\h,\g)\xdif t$ is continuous and bounded as $|\phi(\h,\g)| \leq \norm{\psi}_\infty\norm\h_1\leq\norm{\psi}_\infty$ for each $(\h,\g)\in\HG_p$, $p\geq1$. This confirms the continuity of $\Theta$ because 
	\begin{equation}
		\int_\X \psi(t,\x)\xdif\Theta(\cm) = \int_{\HG_p}\left(\int_0^1\h(t)\psi(t,\g(t))\xdif t\right)\xdif\cm(\h,\g) = \int_{\HG_p}\phi(\h,\g)\xdif\cm(\h,\g)
	\end{equation}
	for each $\cm\in\MC(\HG_p)$, therefore $\Phi=\{\phi\}$ is sufficient. The $p=+\infty$ case is special because for each $t\in\ci01$, the map $\phi(\h,\g)\eqdef \Psi(t,\h,\g)$ is continuous and bounded by $\norm{\psi}_\infty$. This leads to
	\begin{equation}
		\int_{\bar\Omega} \psi(t,\x)\xdif\cmt{t} = \int_{\HG_\infty}\left(\h(t)\psi(t,\g(t))\right)\xdif\cm(\h,\g) = \int_{\HG_\infty}\phi(\h,\g)\xdif\cm(\h,\g),
	\end{equation}
	so we conclude similarly that $\cmt[]{t}$ is continuous in $(\Co(\bar\Omega))'$ as required.
\end{proof}

%% file: lemmas.tex


%% file: energy_results.tex
\section{Results for the structure of \texorpdfstring{$\op{E}$}{E}}\label{app: structure of E}
Here we prove Theorem~\Rref{thm: structure of E} as a result of a sequence of simple lemmas. We start with results for the linear term $\reg$. Throughout this section we fix lower semi-continuous $w\colon\HG\to\ci{0}{+\infty}$ for a complete separable metric space $\HG$ and define $\reg\colon\MC^+(\HG)\to\ci{0}{+\infty}$ by
\begin{equation}
	\forall\cm\in\MC^+(\HG),\qquad \reg(\cm) \eqdef \int_{\HG} w(\hg)\xdif\cm(\hg).
\end{equation}

\subsection{Properties of \texorpdfstring{$\reg$}{W}}
The following lower semi-continuity result is well known (see for instance \cite[Prop. 7.1]{santambrogio_optimal_2015}, or \cite[Prop. 5.1.7]{Ambrosio2008}).
\begin{lmm}[Fatou's lemma for measures]\label{thm: Fatou's lemma}
	If $w\colon \HG\rightarrow \ci0{+\infty}$ is lower semi-continuous, then the functional $\reg$ is lower semi-continuous on $\MC^+(\HG)$ with respect to the narrow topology.
\end{lmm}

Compactness in spaces of measure can be characterised by the following lemmas.
\begin{lmm}[{Prokhorov's theorem on measure spaces, e.g. \cite[Thms. 7.1.7, 8.6.7-8]{Bogachev2007}}]\label{thm: Prokhorov}
	\hfill\\If $\HG$ is a complete separable metric space and $\C A\subset \MC(\HG)$ is a family of Borel measures, then 
	\begin{equation}
		\C A \text{ is relatively compact in the narrow topology \quad if and only if }\quad \C A \text{ is tight and norm-bounded.}
	\end{equation}
\end{lmm}
Because of this lemma, we can give sufficient conditions for the compactness of sub-levelsets of $\reg$.

\begin{lmm}\label{thm: compact W sublevel}
	Let $D$ be a closed subset of $\MC^+(\HG)$ and denote $U_t = \{\cm\in D\st \reg(\cm)\leq t\}$ for $t\in \xR$. If $w\colon \HG \rightarrow \ci0{+\infty}$ is lower semi-continuous, $w$ has compact sub-levelsets, and $U_t$ is bounded in norm for each $t$, then $\reg|_D$ has compact sub-levelsets in the narrow topology of $\MC^+(\HG)$.
\end{lmm}
\begin{proof}
	As $w\geq0$, Lemma~\Rref{thm: Fatou's lemma} shows that $\reg$ is lower semi-continuous, therefore $U_t$ is closed. As $U_t$ is bounded, Lemma~\Rref{thm: Prokhorov} equates compactness with tightness. Finally, \cite[Prop. A.2]{BrediesUnbalanced} (adapted from \cite[Rem. 5.1.5]{Ambrosio2008}) states:
	\begin{quote}
		If $\HG$ is a complete separable metric space and $w$ has compact sub-levelsets, then $U_t$ is tight for each $t\in\xR$.
	\end{quote}
	We conclude that $U_t$ is narrowly compact for each $t\in\xR$.
\end{proof}

We now confirm that the results of this section are applicable to the examples of Sections~\Rref{sec:BB example intro} and \Rref{sec:WFR example}.
\begin{lmm}[Lower semi-continuity and coercivity of optimal transport regularisations]\label{lmma: WFR is coercive and lsc}
	Choose $\alpha,\beta,\delta> 0$, $\HG$ a closed subset of $\HG_\infty$, and let $w:\HG\rightarrow \ho{0}{+\infty}$ be defined by	
	\begin{align}
		\forall (\h,\g)\in \HG,\qquad w(\gamma)&\eqdef \splitln{\int_{\h>0}\left[\alpha + \frac\beta2|\g'|^2 + \frac{\beta\delta^2}2\left(\frac{\h'}{\h}\right)^2\right]\h\xdif t 
		}{\mbox{$\sqrt{\h}, \sqrt{\h}\g \in \ACtwo$}}%
		{+\infty}{\mbox{otherwise.}}
	\end{align}
	Then $w$ is lower semi-continuous and has compact sub-levelsets (in the metric $\metric$).
\end{lmm}
\begin{proof}
	\edit{}{
	It has already been shown that the map $h\delta_\g\mapsto w(\h,\g)$ is lower semi-continuous with compact sub-levelsets in $\enscond{t\mapsto\h(t)\delta_{\g(t)}}{w(\h,\g)<+\infty}$ with the flat metric \cite[Prop. 3.10]{BrediesUnbalanced}. In the proof of Lemma~\Rref{app: metric properties}, we showed that this metric is isometrically equivalent to $\op{d}_\HG$, so $w$ is also lower semi-continuous and coercive with respect to $\op{d}_\HG$.
	}
\end{proof}

Note that the lower semi-continuity and coercivity of the Benamou-Brenier penalty follows immediately if we consider it to be the function
\begin{equation}
	(\h,\g) \mapsto w(\h,\g) + \splitln{0}{\h=1}{+\infty}{\text{else.}}
\end{equation}
The constraint-set is closed, therefore the function is still coercive and lower semi-continuous.

\subsection{Properties of \texorpdfstring{$D$}{D}}
We consider $D$ of the form in \eqref{eq: D def}, that is for some lower semi-continuous $\varphi\colon\HG\to\ho0{+\infty}$,
\begin{equation}\label{app: D def.}
	D \eqdef \left\{\cm\in \MC^+(\HG)\quad\st\quad \int_{\HG}\varphi\xdif\cm\leq1\right\}.
\end{equation}
The only difference between $w$ and $\varphi$ is that $\varphi>0$, so we can re-use many of the results for $\reg$.
\begin{lmm}\label{lmma: extreme point computation}
	If $\varphi\colon\HG\to\ho0{+\infty}$ is lower semi-continuous, then $D$ is closed and
	\begin{equation}\label{eq:extr point D}
		\Ext{D} = \{0\}\cup\enscond{\varphi(\hg)^{-1}\delta_\hg}{\varphi(\hg)<+\infty}.
	\end{equation}
	Furthermore, if $\inf_{\hg\in\HG}\varphi(\hg)\geq\epsilon>0$, then $D$ is bounded.
	Finally, if $\varphi$ has compact sub-levelsets, then $D$ is also compact.
\end{lmm}
\begin{proof}
	Note by Lemma~\Rref{thm: Fatou's lemma} that $D$ is a sub-levelset of the lower semi-continuous function $\cm\mapsto \int_\HG\varphi\xdif\cm$, therefore it is closed. Also, for any $\cm\in\MC^+(\HG)$ with $\norm{\cm}>\epsilon^{-1}$, we have $\int_{\HG}\varphi\xdif\cm >1$. In particular, $\cm\notin D$, so $D$ is bounded. If $\varphi$ has compact sub-levelsets, taking $w=\varphi$ in Lemma~\Rref{thm: compact W sublevel} confirms that $D$ is compact.
	
	It remains to prove \eqref{eq:extr point D}, we begin with the ``$\supset$'' inclusion. By non-negativity, note that for all $\cm, \cm_0, \cm_1 \in D$ and $\lambda\in\oi01$,
	\begin{equation}
		\cm = \lambda\cm_0+(1-\lambda)\cm_1 \qquad\implies \qquad \cm_0,\cm_1\ll \cm, \text{ \ie\ }\op{supp}(\cm_i)\subset\op{supp}(\cm).
	\end{equation}
	In particular, taking $\cm=0$, we obtain $\cm_0=\cm_1=0$, hence $0\in\Ext D$. Now, setting $\cm= \varphi(\hg)^{-1}\delta_\hg$ for some $\hg\in\HG$ with $\varphi(\hg)<+\infty$, we deduce that 
	\begin{align}
		\cm_0=\frac{\alpha}{\varphi(\hg)}\delta_\hg,\ \cm_1=\frac{\beta}{\varphi(\hg)}\delta_\hg,\ \mbox{for some $\alpha,\beta\in\ci01$}.
	\end{align}
	Since $\int_{\HG} \varphi \xdif\cm=1$, we must have $\lambda\alpha+(1-\lambda)\beta=1$ hence $\alpha=\beta=1$, so that $\cm_0=\cm_1=\cm$. As a result, $\varphi(\hg)^{-1}\delta_\hg$ is an extreme point of $D$.
	
	\noindent To show that this inclusion is sharp, we make the following claim.
	\begin{quote}
		\layout{\vspace{-10pt}}
		\begin{clm}
			If $\cm\in D$ and there exists $\hg_0,\hg_1\in\op{supp}(\cm)$ distinct, then $\cm\in\oi{\cm_0}{\cm_1}$ for some $\cm_0,\cm_1\in D$ such that $\hg_0\in\op{supp}(\cm_0)\setminus\op{supp}(\cm_1)$ and $\hg_1\in\op{supp}(\cm_1)\setminus\op{supp}(\cm_0)$.
		\end{clm}
		\begin{proof}[Proof of claim] 
		Let $r = \frac12\metric(\hg_0,\hg_1)$ and set $\HG_{1} = \{\metric(\hg,\hg_1) < r\}$. Define $\cm_1 = \1_{\HG_1}\cm$ and $\cm_0 = \cm-\cm_1$. By the definition of support, for both $i=0,1$ we have $\cm_i\neq 0$, $\hg_i\in\op{supp}(\cm_i)\setminus\op{supp}(\cm_{1-i})$. Also, $\alpha\cm_0+\beta\cm_1\in\MC^+(\HG)$ for all $\alpha,\beta\geq0$, the only challenge is the constraint with $\varphi$.
		\begin{description}
			\item[Case $\bf{\int_{\HG} \varphi\xdif\cm_i = 0}$] In this case, by non-negativity, $\cm_i=0$ so $\op{supp}(\cm_i)=\emptyset$ contradicts the assumption. We conclude that $\int_{\HG} \varphi\xdif\cm_i\in\oi01$ for both $i=0,1$.
			\item[Else] Set $\lambda\eqdef \int_{\HG} \varphi\xdif\cm_0\in\oi01$ and $\int_{\HG} \varphi\xdif\cm_1 = \int_{\HG} \varphi(\xdif\cm-\xdif\cm_0) \leq 1-\lambda\in\oi01$, therefore 
			\begin{equation}
				\cm = \lambda\frac{\cm_0}{\lambda} + (1-\lambda)\frac{\cm_1}{1-\lambda},
			\end{equation}
			which confirms $\cm\in\oi{\cm_0}{\cm_1}$ as required.
		\end{description}
		\end{proof}
	\end{quote}
	
	This claim immediately confirms that all extreme points must have at most one point in their support. Combined with the constraint $\int_\HG\varphi\xdif\cm\leq 1$, we must have 
	\begin{equation}
		\Ext D \subset \{0\} \cup \enscond{\lambda\delta_\hg}{\hg\in\HG,\ \varphi(\hg)<+\infty,\ 0<\lambda\leq \varphi(\hg)^{-1}}.
	\end{equation}
	Finally, if $0<\lambda<\varphi(\hg)^{-1}$, then $\lambda\delta_\hg \in\oi{0}{\varphi(\hg)^{-1}\delta_\hg}$, so $\lambda\delta_{\hg}\notin\Ext{D}$. This confirms that the only extreme points are those found in the first half of the proof.
	
\end{proof}

\subsection{Lower semi-continuity of \texorpdfstring{$\op{E}$}{E}}
Recall that $\op{E}\colon\MC^+(\HG)\to\xR$ is defined by $\op{E}(\cm) = \op{F}(\cm) + \reg(\cm)$. In Theorem~\Rref{thm: structure of E} we assume $\varphi,w\colon\HG\to\ci0{+\infty}$ are lower semi-continuous, therefore Lemmas~\Rref{thm: Fatou's lemma} and \Rref{lmma: extreme point computation} confirm that $\reg$ is lower semi-continuous and $D$ is closed. It remains to show that $\op{F}$ is lower semi-continuous. Lemma~\Rref{lmma: Gamma_infty vs Gamma_1} shows that $\cmt[]{t}$ is narrowly continuous, therefore $\op{F}$ inherits lower semi-continuity from each $\op{F}_j$. From now on we consider any closed subset $\HG\subset\HG_\infty$ with the topology induced by $\metric$.

\subsection{Compactness of sub-levelsets of \texorpdfstring{$\op{E}$}{E}}
We now prove that minimisers of $\op{E}$ exist by showing that $\op{E}$ has compact sub-levelsets. Note in the case $\varphi$ has compact sub-levelsets, then $D$ is already compact (Lemma~\Rref{lmma: extreme point computation}), so the following theorem is not required.
\begin{thrm}\label{thm: compact sub-levelsets}
	Suppose $D\subset\MC^+(\HG)$ is narrowly closed, bounded, and $\op{F}$ is convex lower semi-continuous.
	If $\op{F}$ is bounded from below, $w\colon\HG\to\ci0{+\infty}$ is lower semi-continuous and has compact sub-levelsets, then $\op{E}|_D$ also has compact sub-levelsets.
\end{thrm}
\begin{proof}
	Lemma~\Rref{thm: compact W sublevel} shows that $\reg|_D$ has compact sub-levelsets. Therefore for any $t\in\xR$
	\begin{align}
		\enscond{\cm\in D}{\op{E}(\cm)\leq t} &= \enscond{\cm\in D}{\op{F}(\cm) + \reg(\cm)\leq t}
		\\&\subset \enscond{\cm\in D}{\reg(\cm)\leq t- \inf_{\tilde\cm\in D} \op{F}(\tilde\cm)}.
	\end{align}
	The function $\op{E}$ is lower semi-continuous, so the left-hand side is closed and the right-hand side is compact by assumption. We conclude that the left-hand side is compact for each $t$.
\end{proof}

\subsection{\texorpdfstring{$\op{E}$}{E} has sparse minimisers}
Since the function $\op{F}$ is of the form
\begin{equation}
	\op{F}(\cm) = \sum_{j=0}^T\op{F}_j(A_j\cmt{t_j}) \quad\text{for some}\quad A_j\colon\MC(\HG)\to\xR^\m, \text{ convex }\op{F}_j,
\end{equation}
with $A_j$ as in \eqref{eq: bounded linear kernels}, $\op{F}$ is convex lower semi-continuous. We can therefore use a representer theorem (\eg\ \cite[Cor. 3.8]{Boyer2019}) to demonstrate the sparsity of minimisers.

\begin{lmm}\label{app:lmma: sparse minimisers}
	Suppose $\argmin_{\cm\in D}\op{E}(\cm)\neq \emptyset$ for some $D\subset\MC^+(\HG)$ closed and bounded of the form in \eqref{app: D def.}. Then there exists a minimiser $\cm^*\in D$ such that 
	\begin{align}
		\cm^* &= \sum_{i=1}^{s} a_i \delta_{\gamma_i} \qtext{for some } a_i\geq0,\ \gamma_i \in \HG\ \mbox{with}\ \varphi(\gamma_i)<+\infty,\label{eq:representerthm}
	\end{align}
	where $s\leq \m(T+1)+1$. If in addition $\int \varphi \xdif\cm^*<1$, then $s\leq \m(T+1)$.
\end{lmm}

\begin{proof}
	We reformulate the problem $\min_{\cm \in D}E(\cm)$ as
	\begin{align}
		\min_{\cm \in V} H(\tilde A\sigma) + R(\sigma)
	\end{align}
	where $V$ is the vector space of all $\cm \in \MC(\HG)$ such that $w \in \xLone_{\abs{\cm}}(\HG)$, $R$ is a convex regulariser
	\begin{align}
		R(\sigma)\eqdef W(\cm)+ \chi_D(\cm) \qtext{where} \chi_D(\cm) \eqdef \splitln{0}{\cm\in D}{+\infty}{\text{else,}}
	\end{align}
	with linearly closed level sets. The observation operator $\tilde A\colon V\rightarrow (\xR^\m)^{(T+1)}$, is linear, defined as $\sigma \mapsto (A_j\cmt{t_j})_{0\leq j\leq T}$, and $H\colon (a_0,\ldots,a_T)\mapsto \sum_{j=0}^{T} F_j(a_j)$ is a convex ``fidelity term''.
	
	Set $k\eqdef \m(T+1)$. We claim that there exists a minimiser $\cm^*$ such that 
	\begin{equation}\label{eq: sparse Choquet formulation}
		\cm^* = \sum_{i=1}^{k+1}\theta_i\ep^i,\qtext{with}\sum_{i=1}^{k+1}\theta_i = 1,\text{ and } \forall i,\ \theta_i\geq0,\ \ep^i\in\Ext{D}.
	\end{equation}
	To prove this, fix $t\eqdef R(\cm)$ for some arbitrary $\cm\in \argmin_{\cm\in D}\op{E}(\cm)$ and define
	\begin{equation}\label{eq:reprthm cvx}
		\tilde D \eqdef D\cap\enscond{\cm \in V}{ \IP{w}{\cm}=t}\subset \{R\leq t\}.
	\end{equation}
	We now split the analysis into two cases depending on the value of $t$.
		
	\begin{description}
		\item[Case $t=\inf_V R$]
			In the case $t=\inf_V R=0$, \cite[Cor. 3.8]{Boyer2019} tells us that there exists a minimiser $\cm^*$ which belongs to an elementary face of $\{R\leq t\}$ with dimension at most $k$, therefore also a face of $\tilde D$ with dimension at most $k$. In particular, by Carathéodory's theorem, we can express $\cm^*$ in the form in \eqref{eq: sparse Choquet formulation} but with $\ep^i\in\Ext{\tilde D}$.
		
		Observe that for $t=0$, we can equivalently write $\tilde D$ as
		\begin{equation}
			\tilde D= \enscond{\cm\in\MC^+(\{w=0\})}{\int_\HG\varphi\xdif\cm\leq 1 }\subset D
		\end{equation}
		therefore the extreme points of $\tilde D$ can be computed explicitly. By Lemma~\Rref{lmma: extreme point computation},
		\begin{equation}
			\Ext{\tilde D} = \{0\}\cup\enscond{\varphi(\hg)^{-1}\delta_\hg}{\varphi(\hg)<+\infty,\ w(\hg)=0} \subset \Ext{D},
		\end{equation}
		so \eqref{eq: sparse Choquet formulation} is satisfied.

		\item[Case $t>\inf_V R$] In this case the minimiser $\cm^*$ from \cite[Cor. 3.8]{Boyer2019} belongs to an elementary face of $\{R\leq t\}$ with dimension at most $(k-1)$ (hence also for $\tilde D$). 
		
		Since $\enscond{\cm \in V}{ \IP{w}{\cm}=t}$ is a hyperplane, it is possible to prove (see for instance \cite{dubins_extreme_1962}) that $\cm^*$ belongs to a face of $D$ with dimension at most $k$. The formulation of \eqref{eq: sparse Choquet formulation} is again given by Carathéodory's theorem.
	\end{description}
	We may now deduce \eqref{eq:representerthm} from \eqref{eq: sparse Choquet formulation}. From Lemma~\Rref{lmma: extreme point computation}, we know that the extreme points of $D$ are either $0$ or of the form $\varphi(\hg)^{-1}\delta_\hg$ where $\varphi(\hg)<+\infty$, hence the general form of \eqref{eq:representerthm} where $s\leq k+1$.

	In the special case of $\int \varphi\xdif \cm^* <1$, one of the atoms $\mu^i$ must be $0$. As a result, we may remove it from the sum, so that $s\leq k$.

\end{proof}